\ifx\documentclass\undefined
   \documentstyle[12pt,a4wide,twoside]{article}
\else
   \documentclass[12pt,twoside]{article}
   \usepackage{a4wide}
\fi

\usepackage{amssymb}
\usepackage{graphicx}
\newcommand	{\nc}	{\newcommand}
\nc	{\be}	{\begin{equation}}
\nc	{\ee}	{\end{equation}}
\nc	{\beq}	{\begin{eqnarray}}
\nc	{\eeq}	{\end{eqnarray}}
\nc	{\beqs}	{\begin{eqnarray*}}
\nc	{\eeqs}	{\end{eqnarray*}}

\nc	{\supp}	{{\rm supp}}
\nc	{\diag}	{{\rm diag}}

\nc{\D}		{\displaystyle}
\nc{\SSS}	{\scriptscriptstyle}

\nc{\comment} [1] {}

\newcount\hour \newcount\minute
\hour=\time    \divide\hour by 60
\minute=\hour  \multiply\minute by 60 \advance\minute by -\time \minute=-\minute
\def\twodigits#1{\ifnum #1 < 10{0#1}\else{#1}\fi}

\renewcommand{\comment}[1]{}

\font\Blackbrd=msbm10 scaled 1200       

\newcommand{\jump}[1]   {\mbox{$\big[ \hspace{-0.7mm} \big[ #1
            \big] \hspace{-0.7mm} \big]$} }

\newcommand{\mean}[1]   {\mbox{$\big\{ \hspace{-0.7mm} \big\{ #1
            \big\} \hspace{-0.7mm} \big\}$} }

\newcommand{\sumt}  {\sum_{T\in\caT}}
\newcommand{\sume}  {\sum_{e\in\caE}}

\newcommand{\osc}{\mathop{\rm osc}\nolimits}

\newcommand{\caE}{{\mathcal E}}

\newcommand{\caN}{{\mathcal N}}

\newcommand{\caT}{{\mathcal T}}

\newcommand{\bft}    {{\mathbf t}}

\renewcommand{\le}  {\lesssim}      

\newcommand{\R}     {\mbox{\Blackbrd R}}    
\newcommand{\Poly}  {\mbox{\Blackbrd P}}    

\renewcommand{\div} {{\rm{div} \,}}     

\newtheorem{theorem}{Theorem}[section]
\newtheorem{lemma}[theorem]{Lemma}
\newtheorem{corollary}[theorem]{Corollary}
\newcommand{\bt}{\begin{theorem}}
\newcommand{\et}{\end{theorem}}
\newcommand{\br}{\begin{remark}}
\newcommand{\er}{\end{remark}}
\newcommand{\bc}{\begin{corollary}}
\newcommand{\ec}{\end{corollary}}
\newcommand{\bl}{\begin{lemma}}
\newcommand{\el}{\end{lemma}}
\newcommand{\bp}{\begin{proposition}}
\newcommand{\ep}{\end{proposition}}
\newcommand{\bd}{\begin{definition}}
\newcommand{\ed}{\end{definition}}
\newcommand{\bex}{\begin{example}}
\newcommand{\eex}{\end{example}}

\newtheorem{definition}[theorem]{Definition}
\newtheorem{example}[theorem]{Example}

\newtheorem{remark}[theorem]{Remark}

\newenvironment{proof}{
{\noindent \bf Proof:}}{\quad \hfill \rule{2mm}{2mm}\medskip}

\usepackage{color}
\def\qed{\hbox{\hskip 6pt\vrule width6pt
height7pt depth1pt  \hskip1pt}\bigskip}

\author{Emmanuel
Creus\'e\footnote{Laboratoire Paul Painlev\'e UMR 8524 and INRIA Lille Nord Europe, Universit\'e de Lille 1,
Cit\'e Scientifique, 59655 Villeneuve d'Ascq Cedex, emmanuel.creuse@math.univ-lille1.fr}, Serge
Nicaise\footnote{LAMAV, FR CNRS 2956, Universit\'e de Valenciennes et du Hainaut
Cambr\'esis, , Institut des Sciences et
Techniques de Valenciennes, F-59313 - Valenciennes Cedex 9 France,
serge.nicaise@univ-valenciennes.fr}}

\begin{document}

  \title{A posteriori error estimator based on gradient recovery by averaging for
  convection-diffusion-reaction problems approximated by
   discontinuous Galerkin methods}

  {
\maketitle


 \noindent

{
    \begin{abstract}
 We consider some  (anisotropic and piecewise
constant) convection-diffusion-rea\-ction problems in domains of
$\R^2$,   approximated by a discontinuous Galerkin method with
polynomials of any degree. We propose two  a posteriori error
estimators
 based on  gradient recovery by averaging. It is shown that these estimators give rise
to an upper bound where the constant is explicitly known  up to some
additional terms that guarantees reliability. The lower bound is also
established, one being robust when the convection term (or the
reaction term) becomes dominant. Moreover, the estimator is
asymptotically exact when the recovered gradient is superconvergent.
The reliability and efficiency of the proposed estimators are
confirmed by some numerical tests.
\end{abstract}



 }

    \thispagestyle{plain}

\section{Introduction}
\label{introduction}

The finite element method is   the more popular one that is commonly
used in the numerical realization of different problems appearing in
engineering applications, like the Laplace equation, the Lam\'e
system, the Stokes system, the Maxwell system, etc.... (see
\cite{brenner:94,ciarlet:78,monk98}). More recently discontinuous
Galerkin finite element methods become very attractive since they
present some advantages.  {For example, they allow to perform "$p$
refinement", by locally increasing the polynomial degree of the
approximation if needed. They can moreover use non-conform meshes
allowing hanging-nodes, making the mesh generation  easier for
concrete industrial applications.} We refer to
\cite{arnold:01,cockburn:99}, and the references cited there, for a
good overview on this topic. Adaptive techniques based on a
posteriori error estimators have become indispensable tools for such
methods. For continuous Galerkin finite element methods, there now
exists a vast amount of literature on a posteriori error estimation
for problems in mechanics or electromagnetism and obtaining locally
defined a posteriori error estimates, see for instance the
monographs \cite{AinsworthOden,BS01,Monk03,verfurth:96b}. On the
other hand a similar theory for discontinuous methods is less
developed, let us quote
\cite{becker:03a,creuse:nmpde,creuse:10, houston:03a,houston:07,karaka:03,riviere:03,sun:05}.

Usually upper and lower bounds are proved in order to guarantee the
reliability and the efficiency of the proposed estimator. Most of
the existing approaches  involve constants depending on the shape
regularity of the elements and/or of the jumps in the coefficients;
but these dependences are often not given. Only a small number of
approaches gives rise to estimates with explicit constants, let us
quote
\cite{AinsworthOden,braess:06,fierro:06, LL83,LW04,neittaanmaki:04,nicaise-wohlmuth:05}
for continuous methods.
  For
discontinuous methods, we may cite the recent papers
\cite{ainsworth:06,cochez:07,creuse:10,ern_nic:07,ern:07,lazarov:06}.

 Our goal is here to consider
convection-diffusion-reaction problems   with discontinuous
diffusion coefficients in two-dimensional domains with Dirichlet
boundary conditions approximated by a discontinuous Galerkin method
with polynomials of any degree. Inspired from the paper
\cite{fierro:06}, which treats the case of continuous diffusion
coefficients approximated by a continuous Galerkin method, we
further   derive some a posteriori estimators with an explicit
constant in the upper bound (1 or a similar constant) up to some
additional terms that are usually superconvergent and some
oscillating terms. The approach, called gradient recovery by
averaging \cite{fierro:06} is based on the construction of a
Zienkiewicz/Zhu estimator, namely the difference in an appropriate
norm of $a\nabla_h u_h-G u_h$, where $\nabla_h u_h$ is the broken
gradient of $u_h$ and $G u_h$ is a $H(\div)$-conforming
approximation of this variable. Here special attention has to be
paid due to the assumption that $a$ may be discontinuous. Moreover
the non conforming part of the error is managed using
 a comparison principle from \cite{ern:07} and a
 standard Oswald interpolation operator
 \cite{ainsworth:06,karaka:03}.
 Furthermore using standard inverse inequalities, we show that our
 estimator is locally efficient. Two interests of this approach
 are first the simplicity of the construction of $Gu_h$, and secondly
its superconvergence property (validated by numerical tests). Let us
finally notice that this paper extends our previous one
\cite{creuse:10} in many directions: first we treat
convection-diffusion-reaction problems instead of purely diffusion
ones. Second
 we track the dependence of the constant in the lower
bounds, in particular we show that the natural extension of the
estimator from \cite{creuse:10} yields a lower bound with a constant
that is not robust when the convection and/or the reaction term
becomes dominant. Consequently we introduce another estimator
(adapted from \cite{verfurth:05}) that is robust, keeping
nevertheless an upper bound with an explicit constant.

The schedule of the paper is as follows: We recall in section 2 the
convection-diffusion-reaction problem, its numerical approximation and recall the
comparison principle from \cite{ern:07}. Section 3 is devoted to the
introduction of the first estimator based on gradient averaging and
the proofs of the upper and lower bounds. The upper bound directly
follows from the construction of the estimator and some results from
\cite{fierro:06}, while the lower bound requires the use of some
inverse inequalities and a special construction of $G u_h$. Since
the lower bound is not robust as the   convection and/or the
reaction term becomes dominant we propose an alternative estimator
in section 4 and show its robustness. Finally in section
\ref{sectionnum} some numerical tests are presented  that confirm
the reliability and efficiency of our estimators and the
superconvergence of $Gu_h$ to $a\nabla u$.

Let us finish this introduction with some notation used in the
remainder of the paper: On $D$, the $L^2(D)$-norm will be denoted by
$\|\cdot\|_D$. In  the case $D=\Omega$, we will drop the index
$\Omega$. The usual norm and semi-norm of $H^{s}(D)$ ($s\geq 0$) are
denoted by $\|\cdot\|_{s,D}$ and $|\cdot|_{s,D}$, respectively.
Finally, the notation $a\lesssim b$ and $a\sim b$ means the
existence of positive constants $C_1$ and $C_2$, which are
independent of the mesh size, the coefficients of the operator and
of the quantities $a$ and $b$ under consideration  such that $a\le
C_2b$ and $C_1b\le a\le C_2b$, respectively. In other words, the
constants may depend on the aspect ratio of the mesh,  as well as
the polynomial degree $l$, but they do not depend on the
coefficients of the operator $a, \beta$ and $\mu$ (see below).

\section{The boundary value problem and its discretization}

Let $\Omega$ be a bounded domain of $\R^2$  with a Lipschitz
boundary $\Gamma$ that we suppose to be polygonal. We consider the
following convection-diffusion-reaction problem with   homogeneous
Dirichlet boundary conditions: \be \label{divgradu}
\left\{\begin{array}{rcl}
-\div(a\hbox{~}\nabla u)+\beta\cdot \nabla u+\mu u &=&f \mbox{~in~$\Omega$},\\
u&=&0 \hbox{~on~$\Gamma$}.
\end{array}\right.
\ee

In the sequel, we suppose that $a$ is a symmetric positive definite
matrix which is piecewise constant, namely we assume that there
exists a partition $\cal P$ of $\Omega$ into a finite set of
Lipschitz polygonal domains $\Omega _1,\cdots,\Omega _J$ such that,
on each $\Omega _j$, $a =a _j$ where $a _j$ is a symmetric  positive
definite matrix. Furthermore we assume that $\beta\in H(\div,
\Omega)\cap L^\infty(\Omega)^2$, $\mu\in L^\infty(\Omega)$ and are
such that $\mu-\frac12 \div \beta\geq 0$. If $\mu-\frac12 \div
\beta=0$ on $\omega\subset\Omega$, then we assume that
$\mu=\div\beta=0$ on $\omega$.

The variational formulation of (\ref{divgradu}) involves the
bilinear form
$$ {\it B}(u,v)=\int_\Omega   (a \nabla u \cdot \nabla v+\beta\cdot \nabla u v+\mu u v) , \forall u, v\in H^1_0(\Omega), $$
and the  corresponding energy norm
\[
||| v|||^2=  {\it B}(v,v)=\int_\Omega   (a \nabla v \cdot \nabla
v+(\mu-\frac12 \div\beta)|v|^2), \forall v\in H^1_0(\Omega).  \]

Given $f \in L^2(\Omega)$,  the weak formulation consists in finding
$u \in H^1_0(\Omega)$ such that \be \label{FV} {\it B}(u,v)=(f,v),
\hbox{~}\forall v \in H^1_0(\Omega), \ee where $(f,v)$ means the
$L^2$ inner product in $\Omega$, i.e., $(f,v)=\int_\Omega f v$.
Invoking the positiveness of $a$ and the hypothesis $(\mu-\frac12
\div\beta)\geq 0$, the bilinear form ${\it B}$ is coercive on
$H^1_0(\Omega)$ with respect to the norm $|||\cdot|||$ and this
coerciveness guarantees that problem (\ref{FV}) has a unique
solution by the Lax-Milgram lemma.

\subsection{Discontinuous Galerkin   approximated problem}
\label{s_approx}

Following  \cite{arnold:01,ern:07,karaka:03}, we consider the
following discontinuous Galerkin approximation of our continuous
 problem:
We consider a triangulation $\caT$ made of triangles $T$  whose
edges are denoted by $e$. We
 assume that this triangulation is regular, i.e.,
 for any element $T$, the ratio $\displaystyle \frac{h_{T}}{\rho _{T}}$ is bounded by a constant $\sigma >0$
 independent of $T$ and of the   mesh size $h=\max_{T\in \caT} h_T$,
 where $h_{T}$ is the diameter of $T$ and $\rho _{T}$ the diameter of its largest inscribed ball.
 We further assume that $\caT$ is conforming with the partition $\cal P$ of $\Omega$, i.e.,
 the matrix $a$ being constant on each
  $T \in \caT$, we then denote by
  $a_T$  the value of   $a$ restricted to an element $T$. With each edge  $e$ of an element $T$,
  we associate its length $h_e$ and a unit normal vector $n_{e}$, while $n_T$ stands for the outer unit
  normal vector along $\partial T$.
  $\caE$ (resp. $\caN$) represents the set of edges (resp. vertices) of the triangulation.
In the sequel, we need to distinguish between  edges (resp.
vertices) included into $\Omega$ or into $\Gamma$, in other words,
we set
\beqs \caE_{int} &=&\{e\in\caE: e\subset \Omega\}, \caN_{int} =\{x\in\caN: x\in \Omega\},\\
\caE_{D} &=&\{e\in\caE: e\subset  \Gamma\}, \caN_{D} =\{x\in\caN:
x\in \Gamma\}.\eeqs

Problem (\ref{FV}) is approximated by  the (discontinuous) finite
element  space:
$$X_h=\left\{v_h \in L^2(\Omega)|{v_h}_{|T} \in \Poly_\ell(T), T \in \caT\right\},$$
where $\ell$ is a fixed positive integer.
Later on we also need the continuous counterpart of $X_h$, namely we
introduce
$$
S_h=\left\{v_h \in C(\overline{\Omega})|{v_h}_{|T} \in
\Poly_\ell(T), T \in \caT\right\},$$ as well as
$$
S_{h,1}=\left\{v_h \in C(\overline{\Omega})|{v_h}_{|T} \in
\Poly_1(T), T \in \caT\right\}.$$ We further need
$$X_{h,1}=\left\{v_h \in L^2(\Omega)|{v_h}_{|T} \in \Poly_1(T), T \in \caT \right\}.$$

For our further analysis we need to define some jumps and means
through any   edge  $e\in\caE$ of the triangulation. For
$e\in\caE_{int}$, we denote by $T^+$ and $T^-$ the two elements of
$\caT$ containing $e$. Let $q\in X_h$, we denote by $q^\pm$, the
traces of $q$ taken from $T^\pm$, respectively. Then we define the
mean of $q$ on $e$ by
$$
\mean{q}_\omega=\omega^+(e)q^++\omega^-(e)q^-,
$$
where the nonnegative weights $\omega^\pm(e)$ have to satisfy
$\omega^+(e) +\omega^-(e)=1$. If $\omega^+(e)=\omega^-(e)=1/2$, we
drop the index $\omega$.
 The   jump of $q$ on $e$ is now defined as
follows: \beqs \jump{q}=q^+n_{T^+}+q^-n_{T^-}. \eeqs Remark that the
jump $\jump{q}$ of $q$  is vector-valued.

For $v \in [X_h]^d$, we denote similarly
$$ \mean{v}_\omega=\omega^+(e) v^++\omega^-(e)v^-.$$

For a boundary edge $e$, i. e. $e\in \caE_D$, there exists a unique
element $T^+\in \caT$ such that $e\subset
\partial T^+$. Therefore the mean and jump of $q$ are defined   by
$\mean{q}=q^+$   and $\jump{q}=q^+n_{T^+}$.

For $q \in X_h$,   we define its broken gradient $\nabla_h q$ in
  $\Omega$ by :
$$
 (\nabla_h q)_{|T}=\nabla q_{|T}, \forall T\in \caT.
$$

For $w\in X_h+H^1_0(\Omega)$ and for $\omega\subset \Omega$,  we set
\beqs |||w|||_T^2&=&\int_T (a
\nabla w \cdot \nabla w+(\mu-\frac12 \div\beta)|w|^2), \forall T\in \caT,\\
|||w|||^2_{\omega}&=&\sum_{T\in \caT: T\subset \bar\omega}
|||w|||_T^2,\quad  |||w|||^2=\sumt |||w|||_T^2, \\
\| w \|_{DG,h}&=&|||w||| + \left(\sume h_e^{-1}
\|\jump{w}\|^2_e\right)^{1/2}. \eeqs
Note that $\|\cdot\|_{DG,h}$ is
a norm on  $X_h$.

With these notations,   we define the bilinear form ${\it B}_h(.,.)$
as follows:
\begin{eqnarray} \label{defbh}
       {\it B}_h(u_h,v_h) & = & (a\nabla_h u_h,\nabla_h v_h)
       +((\mu-\div\beta)u_h, v_h)-(u_h,\beta\cdot \nabla_h v_h)
\\&-&
\nonumber \sume\int_e(\theta \mean{ a\nabla_h
v_h}_\omega\cdot\jump{u_h}+\mean{a \nabla_h
u_h}_\omega\cdot\jump{v_h})
\\
\nonumber&+&   \sume \int_e(\gamma_e \jump{u_h}\cdot\jump{v_h}+
\mean{u_h}\beta\cdot \jump{v_h}), \qquad \forall u_h, v_h\in X_h,
\end{eqnarray}
where $\theta$ is   a fixed  real parameter and
 the positive
parameters $\gamma_e$ are chosen appropriately, see below.

The  discontinuous Galerkin approximation of problem
 (\ref{FV}) reads now: Find $u_h \in X_h$,
  such that
\be\label{FVdiscrete} {\it B}_h(u _{h},v _{h}) =(f,v_h), \mbox{~}
\forall v_{h} \in X_h.\ee

In (\ref{defbh}), taking the interior weights $\omega^\pm(e)$ equal
to $1/2$ and setting $\theta=0$,  $\theta=-1$ or $\theta=1$ leads to
the incomplete, nonsymmetric or symmetric interior-penalty
discontinuous Galerkin methods. The stabilization parameters
$\gamma_e$ are chosen in the form (see e.g. \cite{ern:07, ern:07b})
\be\label{choixpoidsstab} \gamma_e=\alpha_e
\frac{\gamma_{a,e}}{h_e}+\gamma_{\beta,e}, \ee where $\alpha_e$ is a
positive parameter, $\gamma_{a,e}$ is a positive parameter that
depends on $a$ and $\gamma_{\beta,e}$ is a positive parameter that
depends on $\beta$ and is zero is $\beta=0$ (the choice
$\gamma_{\beta,e}=\frac{|\beta\cdot n_e|}{2}$ corresponding to an
upwinding scheme). Whenever $\theta\ne -1$, the parameters
$\alpha_e$ have to be large enough to ensure coerciveness of the
bilinear form ${\it B}_h$ on $X_h$ (see, e.g., Lemma 2.1 of
\cite{karaka:03}).

If $\omega$ is a subset of $\Omega$, we denote by $c_{a,\omega}$
(resp. $C_{a,\omega}$) the minimal (resp. maximal) eigenvalue of the
restriction of the matrix $a$ to $\omega$. We further denote by
$c_{\beta,\mu, \omega}=\inf_{x\in \omega}(\mu-\frac12 \div \beta)$.
We recall that from our assumption, we have the implication
\[
c_{\beta,\mu, \omega}=0 \Rightarrow \div \beta=0 \hbox{ and } \mu=0
\hbox{ in } \omega.
\]
For shortness we drop the index $\Omega$ in these notations, namely
we set $c_a=c_{a,\Omega}$, $C_a=C_{a,\Omega}$ and
$c_{\beta,\mu}=c_{\beta,\mu, \Omega}$.

As our approximated scheme is a non conforming one (i.e. the
solution does not belong to $H^1_0(\Omega)$), as usual we need to
take into account the non conforming part of the error. The
possibilities are twofold : either use an appropriate Helmholtz
decomposition of the error (see Lemma 3.2 of \cite{dari:96}, or
Theorem 1 of \cite{ainsworth:06} in 2D or Lemma 2.1 of
\cite{cochez:07}) or use the following comparison principle,  proved
in Theorem 6.3 of \cite{ern:07}. \bl[Abstract upper error bound]
\label{lerrordecompo} Let $u\in H^1_0(\Omega)$ be a solution of
(\ref{FV}) and $u_h \in X_h$ be the solution of (\ref{FVdiscrete}),
then
 \beq\label{helmholtzdec}
&&|||u-u_h|||\leq \inf_{s\in H^1_0(\Omega)} \Big\{|||u_h-s||| \\
&+&\inf_{\bft\in H(\div,\Omega)} \sup_{\varphi\in
H^1_0(\Omega):|||\varphi|||=1 } \Big((f-\div \bft-\beta\cdot \nabla
s-\mu s,\varphi)\nonumber \\
&-&(a\nabla_h u_h+\bft,\nabla \varphi)+((\mu-\frac12\div
\beta)(s-u_h),\varphi)\Big)\Big\}. \nonumber\eeq \el

\section{The  a posteriori error analysis based on gradient recovery by averaging}

Error estimators can be constructed in many different ways as, for
example, using residual type error estimators which measure locally
the jump of the discrete flux \cite{karaka:03}. A different method,
based on equilibrated fluxes, consists in solving local Neumann
boundary value problems \cite{AinsworthOden} or in using
Raviart-Thomas interpolant
\cite{ainsworth:06,cochez:07,ern_nic:07,ern:07}. Here, as an
alternative we introduce a gradient recovery by averaging and define
an error estimator based on a $H(\div)$-conforming approximation of
the broken gradient $a\nabla_h u_h$. In comparison with
\cite{fierro:06}, we admit  discontinuous diffusion coefficients and
use a discontinuous Galerkin method.

Inspired from \cite{fierro:06}  the conforming part of the estimator
$\eta_{CF}$  involves the difference between the broken gradient $a
\nabla_h u_h$ and its smoothed version $Gu_h$, where $Gu_h$ is for
the moment any   element in $X_{h,1}^2$ satisfying
\beq\label{regGuh1} &&Gu_h\in H(\div,\Omega)=\{v\in L^2(\Omega)^2:
\div v\in
L^2(\Omega)\},\\
&& \label{regGuh2} (Gu_h)|_{|\Omega_j}\in H^1(\Omega_j), \forall
j=1,\cdots, J.\eeq
 Hence conforming part of the estimator
$\eta_{CF}$  is defined by \beq \label{defestiRT1}
\eta_{CF}^2=\sum_{T \in \caT} \eta_{CF,T}^2, \eeq where the
indicator $\eta_{CF,T}$ is defined by
$$\eta_{CF,T}=\|a^{-1/2}\left(a\nabla u_h-Gu_h\right) \|_T.$$

For the nonconforming part of the error, we associate  with $u_h$,
its Oswald interpolation operator, namely the unique element $I_{\rm
Os} u_h\in S_h$  defined in the following natural way (see Theorem
2.2 of \cite{karaka:03}): to each node $n$ of the mesh corresponding
to
 the Lagrangian-type degrees of freedom of $S_h$, the value of $I_{\rm Os} u_h$ is the average of the
values of $u_h$ at this node $n$ if it belongs to $\Omega$ (i.e.,
$I_{\rm Os} u_h(n)=\frac{\sum_{n\in T} |T| u_{h|T}(n)}{\sum_{n\in T}
|T|}$) and is zero at this node if it belongs to $\Gamma$. Then the
non conforming indicator $\eta_{NC,T}$ is simply
\[ \eta_{NC,T}=|||I_{\rm Os} u_h-u_h|||_T=
(||a^{1/2} \nabla (I_{\rm Os} u_h-u_h)\|_T^2+||(\mu-\frac12\div
\beta)^{1/2} (I_{\rm Os} u_h-u_h)\|_T^2)^\frac12.
\]
The   non conforming part   of the estimator is then \beq
\label{defestiNC} \eta_{NC}^2=\sum_{T \in \caT} \eta_{NC,T}^2. \eeq

Similarly by keeping in $\eta_{NC}$ only the zero order term, we get
a second non conforming indicator $\eta_{NC2,T}$, which is simply
\[ \eta_{NC2,T}=||(\mu-\frac12\div
\beta)^{1/2} (I_{\rm Os} u_h-u_h)\|_T,
\]
with a global contribution  \beq \label{defestiNC2}
\eta_{NC2}^2=\sum_{T \in \caT} \eta_{NC2,T}^2. \eeq

We obviously notice that
\[
\eta_{NC2,T}\leq \eta_{NC,T}.
\]

Similarly we introduce the estimator corresponding to jumps of
$u_h$: \[\eta_{J}^2=\sume \eta_{J,e}^2, \]
 with
 \[
\eta_{J,e}^2=\left\{
\begin{tabular}{ll}
$ h_e^{-1}\|\jump{u_h}\|_e^2$ &if
$e\in\caE_{int}$,\vspace{2mm}\\
$h_e^{-1}\|u_h\|_e^2$ &if $e\in\caE_{D}$.\end{tabular} \right.
\]

As in \cite{fierro:06}, we introduce some  additional
superconvergent  security parts. In order to   define them properly
we recall that for a node $x\in \caN$, we denote by $\lambda_x$ the
standard hat function (defined as the unique element in $S_{h,1}$
such that $\lambda_x(y)=\delta_{x,y}$ for all $y\in \caN$), let
$\omega_x$ be the patch associated with $x$, which is simply the
support of $\lambda_x$ and let $h_x$ be the diameter of $\omega_x$
(which is equivalent to the diameter $h_K$ of any triangle $K$
included into $\omega_x$). We now denote by $r$ the element residual
$$
r=f+\div(Gu_h)-\beta\cdot \nabla I_{\rm Os} u_h-\mu I_{\rm Os} u_h$$
and for all $x\in \caN$, we set \beqs \bar r_x=&\big(\int_{\omega_x}
\lambda_x\big)^{-1} \int_{\omega_x} r
\lambda_x &\hbox{ if } x\in \caN_{int}=\caN\setminus \caN_D,\\
\bar r_x=&0 &\hbox{ if } x\in \caN_D. \eeqs

We further use a multilevel decomposition of $S_{h,1}$, namely we
suppose  that we start from a coarse grid $\caT_0$ and that the
successive triangulations are obtained by using the bisection
method, see \cite{fierro:06,schmidtsiebert05}. This means that we
obtain a finite sequence of nested triangulations $\caT_\ell$,
$\ell=0, \cdots, L$ such that $\caT_L=\caT$. Denoting by $S_\ell$
the space $$S_{\ell}=\left\{v \in C(\overline{\Omega})|{v}_{|T} \in
\Poly_1(T), T \in \caT_\ell\right\},$$ then we have
\[
S_{\ell}\subset S_{\ell+1} \hbox{ and } S_{h,1}=\cup_{\ell=0}^L
S_{\ell}=S_L.
\]
Furthermore if we denote by $\caN_\ell$ the nodes of the
triangulation $\caT_\ell$, we have
\[
\caN_{\ell}\subset \caN_{\ell+1}.
\]
As usual for all $z\in \caN_{\ell}$ we denote by $\lambda_{\ell z}$
the hat function associated with $z$, namely the unique element in
$S_\ell$ such that
\[\lambda_{\ell z}(z')=\delta_{zz'}  \forall z'\in \caN_{\ell}.
\]

For all $\ell\geq 1$ we finally set
\[
\tilde \caN_{\ell}=(\caN_{\ell}\setminus \caN_{\ell-1})\cup \{z\in
\caN_{\ell-1}\ : \ \lambda_{\ell z}\ne \lambda_{\ell-1 z}\},
\]
and $\tilde \caN_{0}=\caN_{0}$. It should be noticed (see for
instance \cite{fierro:06}) that to each $z\in\tilde \caN_{\ell}$,
the corresponding hat function $\lambda_{\ell z}$  does not belong
to $S_{\ell-1}$.

Now we define $\bar\rho$ and $\bar\gamma$ by \beqs
\bar\rho^2&=&\sum_{x\in \caN} c_{a,\omega_x}^{-1}\rho_x^2,
\\
\bar\gamma^2&=&\sum_{\ell=0}^L\sum_{z\in \tilde\caN_\ell\setminus
\Gamma_D} \gamma_{\ell z}^2, \eeqs where \beqs \rho_x^2&=&
h_x^2\int_{\omega_x} |r-\bar r_x|^2 \lambda_x,
\\
\gamma_{\ell z}&=& |\langle R,\lambda_{\ell z} \rangle|, \eeqs
 $R$ being the residual defined by
\[\langle R,\varphi\rangle=
\int_\Omega  (G u_h\cdot \nabla \varphi+\beta\cdot \nabla I_{\rm Os}
u_h\varphi+\mu I_{\rm Os} u_h\varphi-  f \varphi), \forall \varphi
\in H^1(\Omega).
\]

\subsection{Upper bound} \label{s_upperRT1}

\bt \label{tupperbound}   Let $u\in H^1_0(\Omega)$ be a solution of
problem (\ref{FV}) and let $u_h$ be its discontinuous Galerkin
approximation, i.e.  $u_h \in X_h$ solution of (\ref{FVdiscrete}).
Then there exists $C > 0$ such that \be \label{upperbound} |||
u-u_h|||\leq \eta_{CF}+\eta_{NC}+\eta_{NC2}+C (\bar\rho+\bar
\gamma), \ee and consequently \be \label{upperbound2} \|
u-u_h\|_{DG,h}\leq  \eta_{CF}+\eta_{NC}+\eta_{NC2} +\eta_J +C
(\bar\rho+\bar \gamma). \ee
 \et
\begin{proof}
We first use the estimate (\ref{helmholtzdec}) with $s=I_{\rm Os}
u_h$ and $\bft=-Gu_h$ to get \beqs &&|||u-u_h|||\leq
 \Big\{|||u_h-I_{\rm Os}
u_h||| \\
&+&  \sup_{\varphi\in H^1_0(\Omega):|||\varphi|||=1 } \Big((f+\div
Gu_h-\beta\cdot \nabla I_{\rm Os} u_h-\mu I_{\rm Os}
u_h,\varphi)\nonumber \\
&-&(a\nabla_h u_h-Gu_h,\nabla \varphi)+((\mu-\frac12\div
\beta)(I_{\rm Os} u_h-u_h),\varphi)\Big)\Big\}. \nonumber\eeqs   By
Cauchy-Schwarz's inequality we directly obtain \be ||| u-u_h|||\leq
\eta_{CF}+\eta_{NC}+\eta_{NC2}+ \sup_{\varphi\in
H^1_0(\Omega):|||\varphi|||=1 } |\langle R,\varphi\rangle|. \ee

Using the arguments from Theorem 4.1 of \cite{fierro:06}, we have
\be\label{estimeeresidu} |\langle R,\varphi\rangle|\leq C
(\bar\rho+\bar \gamma) \|a^{1/2}\nabla \varphi\|.\ee

These two estimates lead to the conclusion.
\end{proof}

\br \label{rsuperconvergence}{\rm
 Let us notice that under a
superconvergence property of $||a^{-1/2} (Gu_h-a\nabla u)||$,
$\overline{\rho}$ and $\overline{\gamma}$ will be proved to be
negligible quantities (see Theorem \ref{tlocallowerbd} below), so
that the error is, in this case, asymptotically bounded above by the
estimator without any multiplicative constant. This superconvergence
property is observed in most of practical cases, as for example in
our numerical tests (see section \ref{sectionnum}). Moreover,
theoretical results for different continuous finite element methods
on structured and unstructured meshes have  been established (see
for example \cite{hsww:01,sw:04,Zha01}), but, to our knowledge, not
yet for discontinuous methods for reaction-diffusion-convection
problems on unstructured multi-dimensional meshes. \qed}\er

\subsection{Asymptotic nondeterioration of the smoothed gradient} \label{s_lowerSG}
This subsection is mainly devoted to the estimate of the error
between the smoothed gradient and the exact solution, where we show
that it is bounded by the local error in the $DG$-norm up to an
oscillating term.

We first start with the estimate of some element residuals. On each
element $T$, let us introduce
\[
R_{T}=f+\div(a_T \nabla u_h)-\beta\cdot \nabla  u_h-\mu  u_h.
\]

\bl\label{lestelres} For all $T\in \caT$, one has
\beq\label{estelres1} h_T\|R_T\|_T&\lesssim& \Big(C_{a,T}^{1/2}+h_T
(\frac{\|\beta\|_{\infty,T}}{c_{a,T}^{1/2}}+
\frac{\|\mu\|_{\infty,T}}{c_{\beta,\mu,T}^{1/2}})\Big)
|||u-u_h|||_T+h_T\|f-\Pi_{T}f\|_{T},\\ \label{estelres2}
h_T\|r\|_T&\lesssim&  C_{a,T}^{1/2} \|a^{-1/2}(a\nabla u-Gu_h)\|_T
+h_T \Big(\frac{\|\beta\|_{\infty,T}}{c_{a,T}^{1/2}}+
\frac{\|\mu\|_{\infty,T}}{c_{\beta,\mu,T}^{1/2}}\Big) |||u-I_{\rm
Os}u_h|||_T\\
&+&h_T\|f-\Pi_{T}f\|_{T}, \nonumber\eeq
 where $\Pi_{T}f$  is the $L^2(T)$-orthogonal
projection of $f$
  onto $\Poly_0(T)$.
\el \begin{proof} We start by the first estimate. Its proof is quite
standard, we give it for the sake of completeness. Denoting by
\[
R'_{T}=\Pi_{T}f+\div(a_T \nabla u_h)-\beta\cdot \nabla  u_h-\mu u_h,
\]
we trivially have
\[
R_T=f-\Pi_{T}f+R'_{T}.
\]
Hence it remains to estimate $h_T\|R'_T\|_T$. For that purpose, we
introduce the standard bubble function $b_T$ (see
\cite{verfurth:96b}) and set $w_T=b_TR'_{T}$. Then by a standard
inverse inequality, we have \beqs
\|R'_T\|_T^2&\lesssim& \int_T R'_T w_T\\
&\lesssim& \int_T (\Pi_{T}f+\div(a_T \nabla u_h)-\beta\cdot \nabla
u_h-\mu u_h) w_T\\
&\lesssim& \int_T (\Pi_{T}f-f) w_T +\int_T (f+\div(a_T \nabla
u_h)-\beta\cdot \nabla u_h-\mu u_h) w_T. \nonumber\eeqs Now using
(\ref{divgradu}), we may write  \[
 \|R'_T\|_T^2 \lesssim \int_T (\Pi_{T}f-f) w_T
+\int_T \Big(\div(a_T \nabla (u_h-u))-\beta\cdot \nabla (u_h-u)-\mu
(u_h-u)\Big) w_T\] and by Green's formula, we obtain \be
 \|R'_T\|_T^2\lesssim \int_T (\Pi_{T}f-f) w_T -\int_T \Big(a_T \nabla
(u_h-u)\cdot \nabla w_T+(\beta\cdot \nabla (u_h-u)+\mu (u_h-u))
w_T\Big). \label{res1.1}
 \ee
 By Cauchy-Schwarz's inequality we obtain
 \beqs \|R'_T\|_T^2
&\lesssim&    \|\Pi_{T}f-f\|_T \|w_T\|_T +\|a_T^{1/2} \nabla
(u_h-u)\|_T C_{a,T}^{1/2}\|\nabla w_T\|_T
\\
&+&\left(\|\beta\|_{\infty,T} c_{a,T}^{-1/2}\|a_T^{1/2}\nabla
(u_h-u)\|_T +\|\mu\|_{\infty,T} c_{\beta,\mu,T}^{-1/2}
\|(\mu-\frac12 \div \beta)^{1/2} (u_h-u)\|_T\right) \|w_T\|_T.
 \eeqs
The inverse inequalities \cite{verfurth:96b} \beq   \label{inv1}
\|w_T\|_T \leq    \|R'_T\|_T,  \\
\label{inv2} \|\nabla w_T\|_T
\lesssim  h_T^{-1}  \|R'_T\|_T,
 \eeq
lead to (\ref{estelres1}).

We proceed similarly for the estimate (\ref{estelres2}), namely for
any $T\in \caT$ we set
\[
r'_T=\Pi_Tf+\div(Gu_h)-\beta\cdot \nabla I_{\rm Os} u_h-\mu I_{\rm
Os} u_h,
\]
and therefore
\[
r=f-\Pi_{T}f+r'_{T} \hbox{ on } T.
\]
Hence it remains to estimate $h_T\|r'_T\|_T$. As before we set
$w_T=b_Tr'_{T}$. Then by a standard inverse inequality, we have
\beqs
\|r'_T\|_T^2&\lesssim& \int_T r'_T w_T\\
&\lesssim& \int_T (\Pi_{T}f-f) w_T +\int_T (f+\div(Gu_h)-\beta\cdot
\nabla I_{\rm Os} u_h-\mu I_{\rm Os} u_h) w_T. \eeqs Using
(\ref{divgradu}) we obtain \[ \|r'_T\|_T^2 \lesssim \int_T
(\Pi_{T}f-f) w_T +\int_T \Big(\div(Gu_h-a_T \nabla u)-\beta\cdot
\nabla (I_{\rm Os} u_h-u)-\mu (I_{\rm Os} u_h-u)\Big) w_T,\] and
Green's formula yields
\be\label{res1.2} \|r'_T\|_T^2
 \lesssim \int_T (\Pi_{T}f-f) w_T -\int_T
\Big((Gu_h-a_T \nabla u) \cdot \nabla w_T+(\beta\cdot \nabla (I_{\rm
Os} u_h-u)+\mu (I_{\rm Os} u_h-u)) w_T\Big).
 \ee
  Cauchy-Schwarz's inequality and the inverse inequalities
  (\ref{inv1})-(\ref{inv2})
yield (\ref{estelres2}).
\end{proof}

We go on with the edge residual. \bl\label{lestelresedge} For all
$e\in \caE_{int}$, we set
\[
 \jump{a\nabla_h u_h\cdot n}_e=a_{T^+}\nabla  u_{h|T^+}\cdot n_{T^+}+
a_{T^-}\nabla  u_{h|T^-}\cdot n_{T^-},\] when  $\omega_e=T^+\cup
T^-$. Then
 one has \be\label{tnondeter6}
 h_e^{1/2}
 \|\jump{a\nabla_h u_h\cdot n}_e\|_e   \lesssim
\max_{T\subset\omega_e} \Big(C_{a,T}^{1/2}+h_T
(\frac{\|\beta\|_{\infty,T}}{c_{a,T}^{1/2}}+
\frac{\|\mu\|_{\infty,T}}{c_{\beta,\mu,T}^{1/2}})\Big)
 |||u-u_h|||_{\omega_e}
 +\osc(f,\omega_e), \ee
 where
 $\osc(f,\omega)^2=\sum_{T\subset\omega} h_T^2\|f-\Pi_{T}f\|_{T}^2$.\el
\begin{proof}
 Again the proof is quite standard, we introduce the edge bubble
 function $b_e$ and set $w_e=E(r_e) b_e$, where $E(r_e)$ is an
 extension of $r_e=\jump{a\nabla_h u_h\cdot n}_e$ to $\omega_e$ (see
 \cite{verfurth:96b}).
 By a standard inverse inequality, we have
\[
\|r_e\|_e^2 \lesssim \int_e r_e w_e=\int_e \jump{a\nabla_h
(u_h-u)\cdot n}_e w_e\] By Green's formula, we obtain \beqs
\|r_e\|_e^2
 &\lesssim& \sum_{T\subset\omega_e}\int_T
\left(a_T \nabla (u_h-u)\cdot \nabla w_e +\div(a_T \nabla
(u_h-u))w_e\right).
 \eeqs
Taking into account  (\ref{divgradu}), we get \beqs \|r_e\|_e^2
 &\lesssim& \sum_{T\subset\omega_e}\int_T
\left(a_T \nabla (u_h-u)\cdot \nabla w_e +\div(a_T \nabla u_h)w_e
+(f-\beta\cdot\nabla u-\mu u)w_e \right)\\
&\lesssim& \sum_{T\subset\omega_e}\int_T \Big(a_T \nabla
(u_h-u)\cdot \nabla w_e +(f+\div(a_T \nabla u_h)-\beta\cdot\nabla
u_h-\mu u_h)w_e\\
& &+(\beta\cdot\nabla (u_h-u)+\mu (u_h-u))w_e \Big)
\\
&\lesssim& \sum_{T\subset\omega_e}\int_T \left(a_T \nabla
(u_h-u)\cdot \nabla w_e +R_{T} w_e +(\beta\cdot\nabla (u_h-u)+\mu
(u_h-u))w_e \right).
 \eeqs
 The conclusion follows from Cauchy-Schwarz's inequality,
 (\ref{estelres1}) and the next inverse inequalities
  \beqs
\|w_e\|_T \lesssim   h_e^{1/2} \|r_e\|_e, \forall T\subset\omega_e, \\
\|\nabla w_e\|_T \lesssim  h_e^{-1/2}  \|r_e\|_e, \forall T\subset\omega_e.
 \eeqs
 \end{proof}

To prove our nondeterioration result we also need the next estimate
of the norm of the difference of a vector field $v$ with an
appropriate projection $g(v)$ with the norm of its jumps
 in the interfaces.
 First of all for any vertex $x$ of one $\Omega_j$  and   belonging to more
 than one sub-domain, we introduce the following local notation:
let
 $\Omega_i$, $i=1,\cdots, n$, $n\geq 2$, the sub-domains that have
 $x$ as vertex. We further denote by $n_i$ the unit normal vector
 along the interface $I_i$ between $\Omega_i$ and $\Omega_{i+1}$ (modulo
 $n$ if $x$ is inside the domain $\Omega$) and oriented from $\Omega_i$ and $\Omega_{i+1}$.
 Now we are able to state the following lemma, for a  proof see Lemma
 3.3 of \cite{creuse:10}.

\bl \label{lemmanormequiv} Assume that $x$ is a vertex   of one
$\Omega_j$ and   belonging to more
 than one sub-domain, and use the notations introduced above.
 Then there exists a positive constant $C$
 that depends only on the geometrical situation of the $\Omega_i$'s
 near $x$ such that for all
 $v^{(i)}\in \R^2$,
$i=1,\cdots, n$, there exist vectors $g(v)^{(i)}\in \R^2$,
$i=1,\cdots, n$ satisfying \be\label{condtranshdiv}
(g(v)^{(i+1)}-g(v)^{(i)})\cdot n_i=0, \forall i=1,\cdots, n, \ee and
such that  the following estimate holds \beq \label{normequiv}
\sum_{i=1}^n |v^{(i)}-g(v)^{(i)}| \leq  C   \sum_{i=1}^n |
\jump{v\cdot n}_i |, \eeq where here $|\cdot|$ means the Euclidean
norm and $\jump{v\cdot n}_i$ means the jump of the normal derivative
of $v$ along the interface $I_i$:
$$\jump{v\cdot n}_i=(v^{(i+1)}-v^{(i)})\cdot n_i, \forall i=1,\cdots, n.$$ \el

Using the above lemma, we are now able to prove the asymptotic
nondeterioration of the smoothed gradient if the following choice
for $Gu_h$ is made (we refer to   \cite{creuse:10} for a similar
construction): We
distinguish the following different possibilities for $x\in \caN$.\\
1) First for all vertex $x$ of the mesh (i.e. vertex of at least one
triangle) such that $x$ is inside one $\Omega_j$, we set
\be\label{defZZ} (Gu_h)_{|\Omega_j}(x)=\frac{1}{|\omega_x|}
\sum_{x\in T} |T| a_T \nabla u_{h|T}(x).\ee 2) Second if $x$ belongs
to the boundary of $\Omega$ and to the  boundary of only one
$\Omega_j$ (hence it  does not belong to the  boundary  of another
$\Omega_k$), we define
$(Gu_h)_{|\Omega_j}(x)$ as before. \\
3) If $x$ belongs to an interface between two different sub-domain
$\Omega_j$ and $\Omega_{k}$ but is not a vertex of these
sub-domains, then we denote by $n_{j,k}$ the unit normal vector
pointing from $\Omega_j$ to $\Omega_{k}$ and set $t_{j,k}$ the unit
orthogonal vector of $n_{j,k}$ so that $(n_{j,k}, t_{j,k})$ is a
direct basis of $\R^2$; in that case we set \beq\label{defZZ1}
(Gu_h)_{|\Omega_j}(x)\cdot n_{j,k}&=&(Gu_h)_{|\Omega_k}(x)\cdot
n_{j,k}=\frac{1}{|\omega_x|} \sum_{x\in T} |T| a_T \nabla
u_{h|T}(x)\cdot n_{j,k},\\
(Gu_h)_{|\Omega_j}(x)\cdot t_{j,k}&=&\frac{1}{|\omega_x\cap
\Omega_j|} \sum_{T\subset\Omega_j: x\in T} |T| a_T \nabla
u_{h|T}(x)\cdot t_{j,k},\\
 (Gu_h)_{|\Omega_k}(x)\cdot
t_{j,k}&=&\frac{1}{|\omega_x\cap \Omega_k|} \sum_{T\subset\Omega_k:
x\in T} |T| a_T \nabla u_{h|T}(x)\cdot t_{j,k}.\eeq 4) Finally if
$x$ is a vertex of at least two sub-domains $\Omega_j$, for the sake
of simplicity we suppose that each triangle $T$ having $x$ as vertex
is included into one $\Omega_j$, and  we take \be\label{defZZvertex}
(Gu_h)_{|\Omega_j}(x)=g(v)^{(j)}\forall j\in {\cal J}_x,\ee where
${\cal J}_x=\{j\in\{1,\cdots, J\}: x\in \bar \Omega_j\}$,
$g(v)^{(j)}$ were defined in the previous Lemma \ref{lemmanormequiv}
with here $v$   given by $v=(a_j\nabla u_{h|T}(x))_{j\in {\cal
J}_x}$.

With these choices, we take \be\label{defZZfinal} (Gu_h)_{|\Omega_j}
=\sum_{x\in \caN\cap \bar \Omega_j} (Gu_h)_{|\Omega_j}(x) \lambda_x,
\forall j=1,\cdots, J,\ee where $(Gu_h)_{|\Omega_j}(x)$ was defined
before.

The main point is that by construction $Gu_h$ satisfies the
requirements (\ref{regGuh1}) and (\ref{regGuh2}) but moreover   the
next asymptotic nondeterioration result holds. In order to track the
constants with respect to the diffusion coefficient we first prove
the following technical lemma.

\bl\label{ltechnique} If $\{n,t\}$ is an arbitrary orthonormal basis
of $\R^2$, then for all $T\in \caT$, one has
\be\label{ltechniqueest} |a_T v|\leq \kappa(a_T) (|(a_T v)\cdot n|+2
C_{a,T} | v \cdot t|), \forall v\in \R^2, \ee where we recall that
$\kappa(a_T)$ is the condition number of the matrix $a_T$ defined by
 \[ \kappa(a_T)= C_{a,T} c_{a,T}^{-1}.
 \]\el
\begin{proof}
For shortness we drop the index $T$.

In a first step we take $w\in \R^2$ such that $w \cdot t=0$. Then in
that case we can write
\[
w=(w \cdot n) n,
\]
and therefore
\[
aw=(w \cdot n) an, \quad aw\cdot n=(w \cdot n) (an\cdot n).
\]
By direct calculations, we obtain
\[
|aw|\leq C_a |w \cdot n|, \quad |w \cdot n|\leq c_a^{-1} |aw\cdot
n|,
\]
that leads to \be\label{ltechniqueest2} |a w|\leq \kappa(a) |a
w\cdot n|. \ee

For an arbitrary $v\in \R^2$ we set $w=v-(v\cdot t) t$ that, by
construction, satisfies $w \cdot t=0$. Hence by
(\ref{ltechniqueest2}) we get \beqs |a v|&\leq& |a w|+|v\cdot t|
|at|\\
&\leq& \kappa(a) |a w\cdot n|+C_a|v\cdot t| \eeqs But by the
definition of $w$, one has \[ a w\cdot n=a v\cdot n- (v\cdot t)
at\cdot n,\] and consequently
\[ |a w\cdot n|\leq |a v\cdot n|+C_a|v\cdot t|.\]
The last inequalities leads to (\ref{ltechniqueest}) reminding that
$\kappa(a)\geq 1$.
\end{proof}

Thanks to this lemma we can prove an asymptotic nondeterioration
result  that is similar to Theorem 3.4   of \cite{creuse:10}, where
we here give the explicit dependence of the constants with respect
to the coefficients $a, \beta$ and $\mu$. \bt\label{tnondeter} If
$\ell\leq 2$, then for each element $T \in \caT$ the following
estimate holds \beq\label{lowerbdZZ} &&\|a_T^{-1/2}(Gu_h-a_T\nabla
u)\|_{T} \le c_{a,T}^{-1/2}C_{a,T}\kappa(a_T)\sum_{e\in \caE_{int}:
e\subset \omega_T} h_e^{-1/2}\|\jump{u_h}\|_e+
\\
\nonumber &&+
\left(1+c_{a,T}^{-1/2}\kappa(a_T)\max_{T'\subset\omega_T}
\Big(C_{a,T'}^{1/2}+h_T
(\frac{\|\beta\|_{\infty,T'}}{c_{a,T'}^{1/2}}+
\frac{\|\mu\|_{\infty,T'}}{c_{\beta,\mu,T'}^{1/2}})\Big)\right)
|||u-u_h|||_{\omega_T} \\
&&+c_{a,T}^{-1/2}\kappa(a_T)\osc(f,\omega_T), \nonumber \eeq where
$\omega_T$ denotes the patch consisting of all the triangles of
$\caT$ having a nonempty intersection with $T$.
 \et
\begin{proof}
By the triangle inequality we may write \beqs
\|a_T^{-1/2}(Gu_h-a_T\nabla u)\|_{T}\leq \|a_T^{-1/2}(Gu_h-a_T\nabla
u_h)\|_{T}+\|a_T^{-1/2}(a_T\nabla u_h-a_T\nabla u)\|_{T}
\\
\leq \|a_T^{-1/2}(Gu_h-a_T\nabla u_h)\|_{T}+|||u- u_h|||_{T}. \eeqs
Therefore it remains to estimate the first term of this right-hand
side. For that purpose, since $T\subset \bar \Omega_j$ for a unique
$j\in\{1,\cdots, J\}$, we may write having in mind the assumption $l\leq 2$ :
\[
(Gu_h-a_T\nabla u_h)_{|T}=\sum_{x\in T}
\{(Gu_h)_{|\Omega_j}(x)-a_j\nabla u_{h|T}(x)\}\lambda_x.
\]
As $0\leq \lambda_x\leq 1$, and since the triangulation is regular,
we get \be\label{tnondeter1} \|a_T^{-1/2}(Gu_h-a_T\nabla
u_h)\|_{T}\le c_{a,T}^{-1/2}\sum_{x\in T}
|(Gu_h)_{|\Omega_j}(x)-a_j\nabla u_{h|T}(x) | h_T. \ee We are then
reduced to estimate the factor $|(Gu_h)_{|\Omega_j}(x)-a_j\nabla
u_{h|T}(x) |$ for all nodes $x$ of $T$. For that purpose, we
distinguish four different cases:
\\
1) If $x\in \Omega_j$, then we use an argument similar to the one
from Proposition 4.2 of \cite{fierro:06} adapted to the DG
situation. By the definition of $Gu_h$, we have
\[
(Gu_h)_{|\Omega_j}(x)=\frac{1}{|\omega_x|} \sum_{T'\subset \omega_x}
|T'| a_j \nabla u_{h|T'}(x),
\]
because in this case all $T'\subset \omega_x$  are included into
$\Omega_j$. As a consequence, we obtain \[
(Gu_h)_{|\Omega_j}(x)-a_j\nabla u_{h|T}(x)= \frac{1}{|\omega_x|}
\sum_{T'\subset \omega_x} |T'| a_j (\nabla u_{h|T'}(x)-\nabla
u_{h|T}(x)),
\]
and therefore
\[
\big|(Gu_h)_{|\Omega_j}(x)-a_j\nabla u_{h|T}(x)\big| \le
\sum_{T'\subset \omega_x}  \big|a_j (\nabla u_{h|T'}(x)-\nabla
u_{h|T}(x))\big|.
\]
For each $T'\subset \omega_x$, there exists a path of  triangles of
$\omega_x$, written $T_i, i=0, \cdots, n$ such that \beqs T_0=T,
T_n=T', T_i\ne T_j, \forall i\ne j, \\
T_i\cap T_{i+1} \hbox{ is an common edge } \forall i=1,\cdots, n-1.
\eeqs Hence by the triangle inequality we can estimate
\[
\big|a_j (\nabla u_{h|T'}(x)-\nabla u_{h|T}(x))\big| \leq
\sum_{i=0}^{n-1} \big|a_j (\nabla u_{h|T_{i+1}}(x)-\nabla
u_{h|T_i}(x))\big|.
\]
Now for each term, since $a_j$ is symmetric and positive definite,
using Lemma \ref{ltechnique} above we have \beqs \big|a_j (\nabla
u_{h|T_{i+1}}(x)-\nabla u_{h|T_i}(x))\big| \leq \kappa(a_j)
\big|\{a_j (\nabla u_{h|T_{i+1}}(x)-\nabla
u_{h|T_i}(x))\}\cdot n_i\big| \\
+2 \kappa(a_j) C_{a_j} \big|(\nabla u_{h|T_{i+1}}(x)-\nabla
u_{h|T_i}(x))\cdot t_i\big|, \eeqs
 where $n_i$ is one fixed unit
normal vector along the edge $T_i\cap T_{i+1}$ and $t_i$ is one
fixed unit tangent vector along this edge.
All together we have shown that \[ |(Gu_h)_{|\Omega_j}(x)-a_j\nabla
u_{h|T}(x) | h_T \le h_T\kappa(a_T)\sum_{e \in \caE_{int}: x\in \bar
e} \{|\jump{a\nabla u_h(x)\cdot n}_e|+ C_{a,T}|\jump{\nabla_h
u_h(x)\cdot t}_e|\}.
\]
Using a norm equivalence and an inverse inequality we obtain
\be\label{tnondeter3} |(Gu_h)_{|\Omega_j}(x)-a_j\nabla u_{h|T}(x) |
h_T \le \kappa(a_T) \sum_{e \in \caE_{int}: x\in \bar e} \{
h_e^{1/2}\|\jump{a\nabla_h u_h\cdot n}_e\|_e+C_{a,T}
h_e^{-1/2}\|\jump{u_h}\|_e\}. \ee
\\
2) If the node $x$ belongs to the boundary of $\Omega$ and to the
boundary of a unique $\Omega_j$, since $(Gu_h)_{|\Omega_j}(x)$ is
defined as in the first case, the above arguments lead to
(\ref{tnondeter3}).
\\
3) If $x$ belongs to an interface between two subdomains and is not
a vertex of them, then it is not difficult to show that
\be\label{tnondeter4} |(Gu_h)_{|\Omega_j}(x)-a_j\nabla u_{h|T}(x) |
h_T \le h_T\sum_{e \in \caE_{int}: x\in \bar e} |\jump{a\nabla_h
u_h(x)\cdot n}_e| \ee holds (due to the regularity of the mesh), and
consequently (\ref{tnondeter3}) is still valid.
\\
4) Finally if $x$ is a vertex of different sub-domains $\Omega_j$,
then   Lemma \ref{lemmanormequiv} yields (\ref{tnondeter4}) and
therefore as before we conclude that (\ref{tnondeter3}) holds.

Summarizing the different cases, by (\ref{tnondeter1})  and
(\ref{tnondeter3}), we have \be\label{tnondeter5}
\|a_T^{-1/2}(Gu_h-a_T\nabla u_h)\|_{T}\le
c_{a,T}^{-1/2}\kappa(a_T)\sum_{x\in T} \sum_{e \in \caE_{int}:x\in
\bar e} \{h_e^{1/2}\|\jump{a\nabla_h u_h\cdot n}_e\|_e+
C_{a,T}h_e^{-1/2}\|\jump{u_h}\|_e\}. \ee

 The  first term of this right hand side is   the standard
edge residuals that were estimated in Lemma \ref{lestelresedge},
while the second term in (\ref{tnondeter5}) is part of the DG-norm
and is then kept. Therefore the    estimate  (\ref{tnondeter6}) in
(\ref{tnondeter5}) leads to (\ref{lowerbdZZ}).
\end{proof}

\subsection{Lower bound} \label{s_lowerRT1}

In the spirit of subsection \ref{s_upperRT1} (see Remark
\ref{rsuperconvergence}), we first provide lower bounds where the
error between the gradient of the exact solution and its smoothed
gradient is involved in the right-hand side. In a second step using
the results of the previous section, we give lower bounds with only
the $DG$-norm of the error.

First using the same arguments than in Proposition 4.1 of
\cite{fierro:06}, we have \bt\label{tlocallowerbd} For all
$T\in\caT$, $x\in \caN$ and $\ell \geq 0, z\in \caN_\ell$, we have
\beq\label{locallowerbdineqtr} \eta_{CF,T}&\leq&
\|a_T^{1/2}\nabla(u_h-u)\|_{T}+\|a^{-1/2}(Gu_h-a\nabla u)\|_{T},
\\
\label{locallowerbd2} c_{a,\omega_x}^{-1/2} \rho_x&\le&
\frac{C_{a,\omega_x}^{1/2}}{c_{a,\omega_x}^{1/2}} \|a^{-1/2}(a\nabla
u-Gu_h)\|_{\omega_x} \\
&+& \nonumber
h_x
 \Big(\frac{\|\beta\|_{\infty,\omega_x}}{c_{a,\omega_x}}+
\frac{\|\mu\|_{\infty,\omega_x}}{c_{a,\omega_x}^{1/2}c_{\beta,\mu,\omega_x}^{1/2}}\Big)
|||u-I_{\rm Os}u_h|||_{\omega_x} +c_{a,\omega_x}^{-1/2}\osc(f,\omega_x),\\
\label{locallowerbd3} \gamma_{\ell z}&\le&
C_{a,z}^{1/2}\|a^{-1/2}(Gu_h-a\nabla u)\|_{\omega_{\ell z}} +h_{\ell
z} \Big(\frac{\|\beta\|_{\infty,\omega_{\ell z}}}{c_{a,z}^{1/2}}+
\frac{\|\mu\|_{\infty,\omega_{\ell z}}}{c_{\beta,\mu,z}^{1/2}}\Big)
|||u-I_{\rm Os}u_h|||_{\omega_{\ell z}}.\quad\quad \eeq \et
\begin{proof}
The first estimate is a simple consequence of the triangle
inequality. For the second one, we notice that
\[
\rho_x\leq h_x \|r\|_{\omega_x}.
\]
Hence the estimate (\ref{locallowerbd2}) follows from
(\ref{estelres2}).

For the third estimate by the definition of $\gamma_{\ell z}$,
(\ref{divgradu}) and Green's formula, we have \[ \gamma_{\ell z}
=\int_{\Omega} \left((Gu_h-a\nabla u)\cdot \nabla \lambda_{\ell z}
+(\beta\cdot \nabla (I_{\rm Os}u_h-u)+\mu (I_{\rm
Os}u_h-u))\lambda_{\ell z}\right).
\]
Hence the conclusion follows from Cauchy-Schwarz's inequality
 and using the estimates \[
  \|\lambda_{\ell z}\|\lesssim  h_{\ell z}, \quad
 \|\nabla \lambda_{\ell z}\|\lesssim 1.\]
\end{proof}

For the non conforming part of the estimator, we make use of Lemma
3.5 of \cite{ern:07} (see also Theorem 2.2 of \cite{karaka:03}) to
directly obtain the

\bt\label{tlowerboundNC} Let the assumptions of Theorem
\ref{tupperbound} be satisfied. For each element $T \in \caT$ the
following estimate holds \be\label{upperboundNC} \eta_{NC,T} \le
(C_{a,T}^{1/2}+h_T\|(\mu-\frac12\div\beta)^{1/2}\|_{\infty,T})
\sum_{e\in \caE: e\subset \omega_T} h_e^{-1/2}\|\jump{u_h}\|_e. \ee
\et

\bc\label{clowerboundNC} Let the assumptions of Theorem
\ref{tupperbound} be satisfied. For each element $T \in \caT$ the
following estimate holds \be\label{upperboundNC2} |||u-I_{\rm
Os}u_h|||_{T} \le
(1+C_{a,T}^{1/2}+h_T\|(\mu-\frac12\div\beta)^{1/2}\|_{\infty,T})
\|u-u_h\|_{DG,\omega_T}, \ee where
\[
\|v\|_{DG,\omega}=|||v |||_{\omega}+\left(\sum_{e\in \caE: e\subset
\omega} h_e^{-1}\|\jump{v}\|_e^2\right)^{1/2}.
\]
\ec

\br\label{rsuperconv1}{\rm From the estimates (\ref{locallowerbd2}),
(\ref{locallowerbd3}) and (\ref{upperboundNC2}), we can say that the
quantities $c_{a,\omega_x}^{-1/2} \rho_x$, $\bar \rho$ and
$\gamma_{L z}$ are superconvergent if $\|a^{-1/2}(a\nabla u-Gu_h)\|$
and $\osc(f,\omega_x)$ are. Obviously this superconvergence
properties are lost when the diffusion matrix $a$ tends to $0$,
i.e., it is not robust with respect to convection dominance. A
similar phenomenon also holds for $\bar \gamma$, see Remark
\ref{rsuperconv2}. \qed} \er

A direct consequence of these three Theorems is the next local lower
bound:

\bt\label{tlowerbound} Let the assumptions of Theorems
\ref{tupperbound} and \ref{tnondeter} be satisfied. For each element
$T \in \caT$ the following estimate holds \beqs
&&\eta_{CF,T}+\eta_{NC;T}+\eta_{J,T}+\sum_{x\in
T}(c_{a,\omega_x}^{-1/2} \rho_x+ \gamma_x) \le \kappa_{1,T}
\|u-u_h\|_{DG,\tilde \omega_T}
\\
 && \hspace{2cm}+ C_{a,\omega_T}^{1/2}(1+c_{a,\omega_T}^{-1/2}) \|a^{-1/2}(a\nabla
u-Gu_h)\|_{\omega_T}+
 \sum_{x\in
T}c_{a,\omega_x}^{-1/2}  \osc(f, \omega_x), \eeqs
 where $\gamma_x=
\gamma_{L x}$ recalling that  $L$ is such that $\caN_L=\caN$,
$\tilde \omega_T$ is the larger patch defined by $\tilde
\omega_T=\displaystyle \cup_{T'\subset \omega_T} {\omega_{T'}}$ and
\[
\kappa_{1,T}=\big(1+C_{a,\omega_T}^{1/2}+h_T\|(\mu-\frac12\div\beta)^{1/2}\|_{\infty,\omega_T}\big)
\Big\{1+h_T(1+c_{a,\omega_T}^{-1/2})
 \Big(\frac{\|\beta\|_{\infty,\omega_T}}{c_{a,\omega_T}^{1/2}}+
\frac{\|\mu\|_{\infty,\omega_T}}{c_{\beta,\mu,\omega_T}^{1/2}})\Big)
\Big\}.
\]
\et

Using the nondeterioration result from the previous subsection, we
get a lower bound with the $DG$-norm of the error at any event.

\bc\label{ctlowerbound} Let the assumptions of Theorems
\ref{tupperbound} and \ref{tnondeter} be satisfied. For each element
$T \in \caT$ the following estimate holds \[
\eta_{CF,T}+\eta_{NC,T}+\eta_{J,T}+\sum_{x\in
T}(c_{a,\omega_x}^{-1/2} \rho_x+ \gamma_x) \le \tilde \kappa_{1,T}
\|u-u_h\|_{DG,\tilde \omega_T}+ c_{a,\tilde\omega_T}^{-1/2}
\kappa(a_{\omega_T})^2 \osc(f,\tilde\omega_T),
\] where
\beqs &&\tilde
\kappa_{1,T}=\kappa_{1,T}+C_{a,\omega_T}^{1/2}(1+c_{a,\omega_T}^{-1/2})
\Big(1+C_{a,\omega_T}^{1/2}\kappa(a_{\omega_T})^2 \\
&&\hspace{2cm}+ c_{a,\omega_T}^{-1/2}\kappa(a_{\omega_T})
\big(C_{a,\tilde \omega_T}^{1/2}+h_T (\frac{\|\beta\|_{\infty,\tilde
\omega_T}}{c_{a,\tilde \omega_T}^{1/2}}+
\frac{\|\mu\|_{\infty,\tilde \omega_T}}{c_{\beta,\mu,\tilde
\omega_T}^{1/2}})\big)\Big). \eeqs \ec

Combining the arguments of  Proposition 4.3 of \cite{fierro:06} and
the above ones, one can obtain a global lower bound:
\bt\label{tlowerbound2} Let the assumptions of Theorems
\ref{tupperbound} and \ref{tnondeter} be satisfied. Then the
following global lower bound   holds \[
\eta_{CF}+\eta_{NC}+\eta_{J}+\bar\rho+\bar\gamma \le
\kappa_2\|u-u_h\|_{DG,h}+\osc_1(f,\Omega), \] where \beqs
\kappa_1&=& (1+C_{a}^{1/2}+ \max_{x\in \caN}
\frac{C_{a,\omega_x}^{1/2}}{c_{a,\omega_x}^{1/2}}),\\
\kappa_2&=& \kappa_3
\\
&+& \kappa_1\Big\{1+\max_T \left(c_{a,T}^{-1/2}C_{a,T}\kappa(a_T)+
c_{a,T}^{-1/2}\kappa(a_T)\max_{T'\subset\omega_T}
\Big(C_{a,T'}^{1/2}+h_T
(\frac{\|\beta\|_{\infty,T'}}{c_{a,T'}^{1/2}}+
\frac{\|\mu\|_{\infty,T'}}{c_{\beta,\mu,T'}^{1/2}})\Big)\right)\Big\},
\\
\kappa_3&=&(1+C_a^{1/2})\Big(1+\max_{x\in \caN} h_x
 (\frac{\|\beta\|_{\infty,\omega_x}}{c_{a,\omega_x}}+
\frac{\|\mu\|_{\infty,\omega_x}}{c_{a,\omega_x}^{1/2}c_{\beta,\mu,\omega_x}^{1/2}})
+c_{\beta,\mu}^{-1/2}(\|\beta\|_{\infty}+\|\div \beta\|_{\infty}+
\|\mu\|_{\infty}) \Big),\\
\osc_1(f,\Omega)&=& \left(\sum_{x\in \caN}c_{a,\omega_x}^{-1}
\osc(f, \omega_x)^2\right)^{1/2} + \kappa_1\left(\sum_{T\in
\caT}c_{a,T}^{-1}\kappa(a_T)^2 \osc(f, \omega_T)^2\right)^{1/2},
\eeqs \et
\begin{proof} As in Proposition 4.3 of \cite{fierro:06}, we have
\[
\bar\gamma^2=\langle R, \chi\rangle, \] with $\chi\in H^1_0(\Omega)$
defined in the proof of Proposition 4.3 of \cite{fierro:06}. As
before using
  (\ref{FV}), we have
\be\label{fierro}
  \langle R, \chi\rangle
=\int_{\Omega} \left((Gu_h-a\nabla u)\cdot \nabla \chi +(\beta\cdot
\nabla (I_{\rm Os}u_h-u)+\mu (I_{\rm Os}u_h-u))\chi\right)\ee and by
Green's formula, we get \be\label{fierro2}
  \langle R, \chi\rangle
=\int_{\Omega} \left((Gu_h-a\nabla u)\cdot \nabla \chi + (I_{\rm
Os}u_h-u) \div(\beta\chi)+\mu (I_{\rm Os}u_h-u)\chi\right).\ee By
Cauchy-Schwarz's inequality we get \[ \langle R, \chi\rangle
 \lesssim
C_{a}^{1/2}\|a^{-1/2}(Gu_h-a\nabla u)\| \|\nabla \chi\| +
\Big(\|\beta\|_{\infty}+\|\div \beta\|_{\infty}+
\|\mu\|_{\infty}\Big) ||u-I_{\rm Os}u_h|| \| \chi\|_{1,\Omega}. \]
Hence by Poincar\'e's inequality, we deduce that
\[
\bar\gamma^2\le \left(C_{a}^{1/2}\|a^{-1/2}(Gu_h-a\nabla u)\|
 + \Big(\|\beta\|_{\infty}+\|\div \beta\|_{\infty}+
\|\mu\|_{\infty}\Big) ||u-I_{\rm Os}u_h||\right)
 \|\nabla \chi\|. \] Since Proposition 4.3 of
\cite{fierro:06} shows that $\|\nabla \chi\|\le \bar\gamma$, we
conclude that \be\label{fierro3} \bar\gamma\le
C_{a}^{1/2}\|a^{-1/2}(Gu_h-a\nabla u)\|
 + \Big(\|\beta\|_{\infty}+\|\div \beta\|_{\infty}+
\|\mu\|_{\infty}\Big) ||u-I_{\rm Os}u_h||. \ee

This estimate and  Theorems \ref{tnondeter}, \ref{tlocallowerbd} and
\ref{tlowerboundNC} lead to the conclusion.
\end{proof}

\br\label{rsuperconv2}{\rm From the estimate (\ref{fierro3})  we can
say that the quantity $\bar \gamma$ is superconvergent if
$\|a^{-1/2}(a\nabla u-Gu_h)\|$ and $||u-I_{\rm Os}u_h||$ are
superconvergent. For the second term $||u-I_{\rm Os}u_h||$, using a
triangle inequality, we have
\[
||u-I_{\rm Os}u_h||\leq  ||u- u_h||+||u_h-I_{\rm Os}u_h||,
\]
and Lemma 3.5 of    \cite{ern:07} yields (see Theorem
\ref{tlowerboundNC})
\[
||u_h-I_{\rm Os}u_h||\le h\|(\mu-\frac12\div\beta)^{1/2}\|_{\infty}
\|u-u_h\|_{DG},
\]
Hence $||u_h-I_{\rm Os}u_h||$ is a superconvergent quantity. On the
other hand using an Aubin-Nitsche trick (see for instance
\cite{arnold:01}), if $\theta=1$, then we have \[ \|u-u_h\|\leq C
h^{\sigma} \|u-u_h\|_{DG}, \] for some $\sigma>0$ and $C$ depends on
the matrix $a$, $\beta$ and $\mu$ (and may blow up as $a\to 0$).
Hence, if $\theta=1$, $\|u-u_h\|$ is  is a superconvergent quantity
and therefore  $||u-I_{\rm Os}u_h||$ as well. As in Remark
\ref{rsuperconv1}, this property is not   valid in the convective
dominant case. \qed} \er

According to Theorem \ref{tlowerbound2} we see that our lower bound
is not robust with respect to  convection/reaction dominance. Indeed
two factors (namely the factors $c_{a,\omega_x}^{-1}\|\beta\|_{\infty,\omega_x}$
and
 $c_{a,\omega_x}^{-1/2}\|\mu\|_{\infty,\omega_x}$ in the definition of $\kappa_3$) blow up
as $a$ goes to zero. The two factors come from two different terms
in the proof of the lower bound. Indeed the first one appears when
one wants to estimate a convective derivative (see Lemma
\ref{lestelres} for instance) and is usually eliminated by adding to
the norm of the error an additional dual norm (see
\cite{schotzau:07,verfurth:05}), this technique will be adapted to
our problem in the next section. We further see that the same
factors appear in  the estimate of the error between $Gu_h$ and
$a\nabla u$, hence a simple  idea   to eliminate these factors is to
add the term $||a^{-1/2}(Gu_h-a\nabla u)||$ to the error (even if in
practice this error is quite often superconvergent). These two ideas
are developed below.

Note also that in Theorem \ref{tlowerbound2} the lower bound is not
robust with respect to the local anisotropy of the diffusion matrix
(via the constants $\kappa_1$ and $\kappa_2$ that blow up if the
condition number $\kappa(a_j)$ blows up for at least one $j$). But
this is a normal phenomenon that also appears in former works, see
for instance Theorem 3.2 of \cite{ern:07}, Theorem 4.4 of
\cite{Vohralik:07} or Theorem 4.2 of \cite{Vohralik:08}.

\section{A  robust a posteriori estimate}

According to the previous papers \cite{schotzau:07,verfurth:05} we
need to use an appropriate   dual norm in order to estimate
convective derivatives. The price to pay is that the constant in the
upper bound will be no more 1, but remains nevertheless explicit.
Note further that theoretically we lose the robust superconvergence
property of $\bar \rho$ (or a modification of it, see below) and of
$\bar \gamma$.

We  start with the following definition: for $h\in L^2(\Omega)$, let
us denote by $|||h|||_*$ its norm as an element of the dual of
$H^1_{0}(\Omega)$ (equipped with the norm $|||\cdot|||$), namely
\[
|||h|||_*=\sup_{v\in H^1_{0}(\Omega)\setminus\{0\}}
\frac{\int_\Omega hv}{|||v|||}.
\]
Then the proof of Lemma 3.1 of \cite{verfurth:05} yields the
following result. \bl\label{lderconv} If
\[
M=\max\{1, \sup_{\Omega}\frac{|\mu|}{(\mu-\frac12\div\beta)}\},
\]
then we have
\[
\left|\int_\Omega (a\nabla v\cdot \nabla w+\mu vw)\right|\leq M\,
|||v||| \, |||w|||, \forall v,w\in H^1_{0}(\Omega),
\]
and
\[
|||w|||+|||\beta\cdot \nabla w|||_*\leq (2+M) \sup_{v\in
H^1_{0}(\Omega)\setminus\{0\}} \frac{B(w,v)}{|||v|||}, \forall w\in
H^1_{0}(\Omega).
\]
\el

For robustness reasons, we further introduce \beqs
\alpha_x&=&\min\{c_{\beta,\mu}^{-1/2}, c_{a,\omega_x}^{-1/2} h_x\},
\forall x\in \caN,\\ \alpha_T&=&\min\{c_{\beta,\mu}^{-1/2},
c_{a_T}^{-1/2} h_T\}, \forall T\in \caT, \eeqs and
\[
\tilde \rho^2 = \sum_{x\in \caN} \alpha_x^2 \int_{\omega_x} |r-\bar
r_x|^2 \lambda_x.
\]

\bt \label{tupperboundrob2}   Let $u\in H^1_0(\Omega)$ be a solution
of problem (\ref{FV}) and let $u_h$ be its discontinuous Galerkin
approximation, i.e.  $u_h \in X_h$ solution of (\ref{FVdiscrete}).
Then there exists $C > 0$ such that \be \label{upperboundrob2} |||
u-u_h|||+|||\beta\cdot \nabla (u-I_{\rm Os}u_h)|||_*\leq
\eta_{NC2}+(3+M)\left(\eta_{CF}+\eta_{NC}+ C (\tilde\rho+\bar
\gamma)\right), \ee and consequently \be  \|
u-u_h\|_{DG,h}+|||\beta\cdot \nabla (u-I_{\rm Os}u_h)|||_*\leq
 \eta_J+\eta_{NC2}+(3+M)\left(\eta_{CF}+\eta_{NC} +C
 (\tilde\rho+\bar
\gamma)\right). \ee
 \et
\begin{proof}
Applying Lemma \ref{lderconv} to $w=u-I_{\rm Os}u_h$, we see that
\[
|||\beta\cdot \nabla (u-I_{\rm Os}u_h)|||_*\leq (2+M) \sup_{v\in
H^1_{0}(\Omega)\setminus\{0\}} \frac{B(u-I_{\rm Os}u_h,v)}{|||v|||}.
\]
But according to (\ref{FV}) and the definition of $R$ we have
\[
B(u-I_{\rm Os}u_h,v)=(a^{-1/2}(Gu_h-a\nabla_h u_h), a^{-1/2} \nabla
v)+(a^{1/2}(\nabla_h u_h-\nabla I_{\rm Os}u_h), a^{-1/2} \nabla
v)-\langle R,v\rangle.
\]
Hence by Cauchy-Schwarz's inequality we deduce that
\[
B(u-I_{\rm Os}u_h,v)\leq (\eta_{CF}+\eta_{NC})|||v||| +|\langle
R,v\rangle|.
\]
Similarly to Theorem 4.1 of \cite{fierro:06}, we can show that
\be\label{estimeeresidu2} |\langle R,\varphi\rangle|\leq C
(\tilde\rho+\bar \gamma) ||| \varphi|||,\ee due to the easily
checked estimate
\[
\|\varphi-\Big(\int_{\omega_x}\lambda_x\Big)^{-1}
\int_{\omega_x}\varphi \lambda_x\|\lesssim \alpha_x |||
\varphi|||_{\omega_x}.
\]
The estimate (\ref{estimeeresidu2}) then yields
\[
B(u-I_{\rm Os}u_h,v)\leq (\eta_{CF}+\eta_{NC}+C (\tilde\rho+\bar
\gamma))|||v|||,
\]
and therefore \be \label{upperbound3} |||\beta\cdot \nabla (u-I_{\rm
Os}u_h)|||_*\leq (2+M)(\eta_{CF}+\eta_{NC}+C (\tilde\rho+\bar
\gamma)). \ee

In the same manner, the estimate (\ref{estimeeresidu2}) allows to
show that
\[ ||| u-u_h|||\leq \eta_{CF}+\eta_{NC}+\eta_{NC2}+C
(\tilde\rho+\bar \gamma). \]
 This estimate and (\ref{upperbound3}) lead to
(\ref{upperboundrob2}).
\end{proof}

As suggested before we now add the error between the gradient of $u$
and its smoothed gradient in order to have a full robust lower
bound. \bc \label{ctupperboundrob2} Under the assumptions of Theorem
\ref{tupperboundrob2},  there exists $C > 0$ such that \beqs
\label{upperboundrob2c} ||| u-u_h|||+\|a^{-1/2}(Gu_h-a\nabla
u)\|+|||\beta\cdot \nabla (u-I_{\rm Os}u_h)|||_*\leq
\eta_{CF}+2\eta_{NC2}\\
+2(3+M)\left(\eta_{CF}+\eta_{NC}+ C (\tilde\rho+\bar\gamma)\right),
\nonumber\eeqs and consequently \beqs  \|
u-u_h\|_{DG,h}+\|a^{-1/2}(Gu_h-a\nabla u)\|+|||\beta\cdot \nabla
(u-I_{\rm Os}u_h)|||_*\leq \eta_J+\eta_{CF}+2\eta_{NC2}
\\
+2(3+M)\left(\eta_{CF}+\eta_{NC}+ C (\tilde\rho+\bar \gamma)\right).
\eeqs
 \ec

The lower bound requires to revisit some results from subsection
\ref{s_lowerRT1}. For that purpose and only for the sake of
simplicity we require that the diffusion matrix $a$ is constant on
the whole domain. Let us now revisit Lemma \ref{lestelres}:

\bl\label{lestelresrob} Assume that the diffusion matrix $a$ is
constant. Then the next estimate holds \beq \label{estelres2rob}
&&\tilde\rho \lesssim \|a^{-1/2}(Gu_h-a\nabla u)\|+|||\beta\cdot
\nabla (u-I_{\rm Os}u_h)|||_*
+\frac{\|\mu\|_\infty}{c_{\beta,\mu}^{1/2}}  \|u-I_{\rm
Os}u_h\|+\xi, \eeq where we have set
\[
\xi^2=\sumt \alpha_T^2\|f-\pi_Tf\|_T^2.
\]
 \el \begin{proof}
First we notice that
\[
\tilde\rho^2\lesssim \sum_{x\in \caN} \alpha_x^2
\|r\|^2_{\omega_x}\lesssim \sum_{T\in \caT} \alpha_T^2
  \|r\|^2_{T}.
  \]
Hence it remains to estimate
\[
\sumt \alpha_T^2 \|r'_T\|_T^2.
\]
For that purpose, we start as in the proof of Lemma \ref{lestelres}.
Namely for all elements $T$, the estimate (\ref{res1.2}) holds.
Hence multiplying this estimate by $\alpha_T^2$, summing on $T$ and
setting $w=\sumt \alpha_T^2 w_T$ (with $w_T$ defined in the proof of
Lemma \ref{lestelres}), we have \beqs \sumt \alpha_T^2 \|r'_T\|_T^2
&\lesssim&    \sumt \alpha_T^2 \int_T (\Pi_{T}f-f) w_T
\\
&-&\int_\Omega \Big(Gu_h-a\nabla u)\cdot \nabla w+(\beta\cdot \nabla
(I_{\rm Os} u_h-u) w+\mu (I_{\rm Os} u_h-u)w\Big). \nonumber
 \eeqs
 By Cauchy-Schwarz's inequality and the definition of the norm
 $|||\cdot |||_*$, we get
\beqs \sumt \alpha_T^2 \|r'_T\|_T^2 &\lesssim&    \xi^2 (\sumt
\alpha_T^2 \|w_T\|_T^2)+ (\|a^{-1/2}(Gu_h-a\nabla u)\|
\\
&+&|||\beta\cdot \nabla (u-I_{\rm
Os}u_h)|||_*+\frac{\|\mu\|_\infty}{c_{\beta,\mu}^{1/2}} \|u-I_{\rm
Os}u_h\|)|||w|||. \nonumber
 \eeqs
By the estimate (4.15) of \cite{verfurth:05}, we see that
\[
|||w|||^2\lesssim \sumt \alpha_T^2 \|r'_T\|_T^2,
\]
while we directly check that
\[
\sumt \alpha_T^2 \|w_T\|_T^2 \leq \sumt \alpha_T^2 \|r'_T\|_T^2.
\]
These three estimates yield
\[
\Big(\sumt \alpha_T^2\|r'_T\|_T^2\Big)^{1/2}\lesssim
\|a^{-1/2}(Gu_h-a\nabla u)\|+|||\beta\cdot \nabla (u-I_{\rm
Os}u_h)|||_* +\frac{\|\mu\|_\infty}{c_{\beta,\mu}^{1/2}}\|u-I_{\rm
Os}u_h\|+\xi,
\]
and we obtain (\ref{estelres2rob}) by the definition of $r'_T$.
\end{proof}

As in Theorem \ref{tlowerbound2}, the above results allow to obtain
the following robust global lower bound.

\bt\label{tlowerbound2rob} Let the assumptions of Theorems
\ref{tupperbound} and \ref{tnondeter} be satisfied. Assume that the
diffusion matrix $a$ is constant. Then the following global lowerÄ
bound holds \beqs
&&\eta_{CF}+\eta_{NC}+\eta_{J}+\tilde\rho+\bar\gamma  \le
\kappa_4|||u-u_h|||+\kappa_5\eta_J \\
&&\hspace{2cm}+ \kappa_6 |||\beta\cdot \nabla (u-I_{\rm Os}u_h)|||_*
+
 \|a^{-1/2}(Gu_h-a\nabla u)\|+\kappa(a) (1+C_a^{1/2})\xi, \eeqs where
\beqs
\kappa_4&=&M  (1+c_{\beta,\mu}^{1/2}),\\
\kappa_5&=&1+\kappa_4(C_a^{1/2}+h\|(\mu-\frac12\div \beta)^{1/2}\|_\infty),\\
\kappa_6&=&1+C_a^{1/2}+\|(\mu-\frac12\div \beta)^{1/2}\|_\infty.
\eeqs
 \et
\begin{proof} By using the identity (\ref{fierro}) of the proof of  Theorem \ref{tlowerbound2}, we
get \beqs \bar\gamma^2&=&\langle R, \chi\rangle\\
 &\lesssim& C_{a}^{1/2}\|a^{-1/2}(Gu_h-a\nabla u)\| \|\nabla \chi\| +
 \|\mu\|_{\infty}  \|u-I_{\rm Os}u_h\| \|
\chi\|+ |||\beta\cdot \nabla (u-I_{\rm Os}u_h)|||_* |||\chi|||.
\eeqs But Poincar\'e's inequality yield
\[
|||\chi|||^2\lesssim (C_a+ \|\mu-\frac12 \div \beta\|_\infty)
\|\nabla \chi\|^2.
\]
These two estimates and the fact   $\|\nabla \chi\|\le \bar\gamma$
(see Proposition 4.3 of \cite{fierro:06}) yield \beqs \bar\gamma
 &\lesssim& C_{a}^{1/2}\|a^{-1/2}(Gu_h-a\nabla u)\|  +
 \|\mu\|_{\infty} \|u-I_{\rm
Os}u_h\| \\
&+& (C_a^{1/2}+ \|(\mu-\frac12 \div \beta)^{1/2}\|_\infty)
|||\beta\cdot \nabla (u-I_{\rm Os}u_h)|||_*. \eeqs

This estimate,  (\ref{locallowerbdineqtr}),  (\ref{estelres2rob})
and (\ref{upperboundNC})  lead to the conclusion.
\end{proof}

\br {\rm The factor $M$ in the lower bound is also present in
Theorem 4.1 of \cite{verfurth:05}. \qed} \er

\section{Numerical results}
\label{sectionnum}
Here we illustrate and validate our theoretical results by some computational experiments.
\subsection{The homogeneous case}
\label{regular}
We consider the domain $\Omega=\{ 0<x,y<1\}$, the reaction coefficient $\mu=1$, the velocity field $\beta=(1,0)^\top$ and the isotropic homogeneous diffusion tensor $a=\varepsilon \, \mathcal{I}$ where $\mathcal{I} $ is the identity matrix. Here, as in \cite{ern:07}, we take $\varepsilon=1E-02$ and $\varepsilon=1E-04$. The source term $f$ is chosen accordingly so that $u=\frac12 x(x-1)y(y-1)(1-tanh(10-20x))$ is the exact solution (see Figure \ref{exact_solution}). \\

\begin{figure}[htp]
\begin{center}
\includegraphics[scale=0.50]{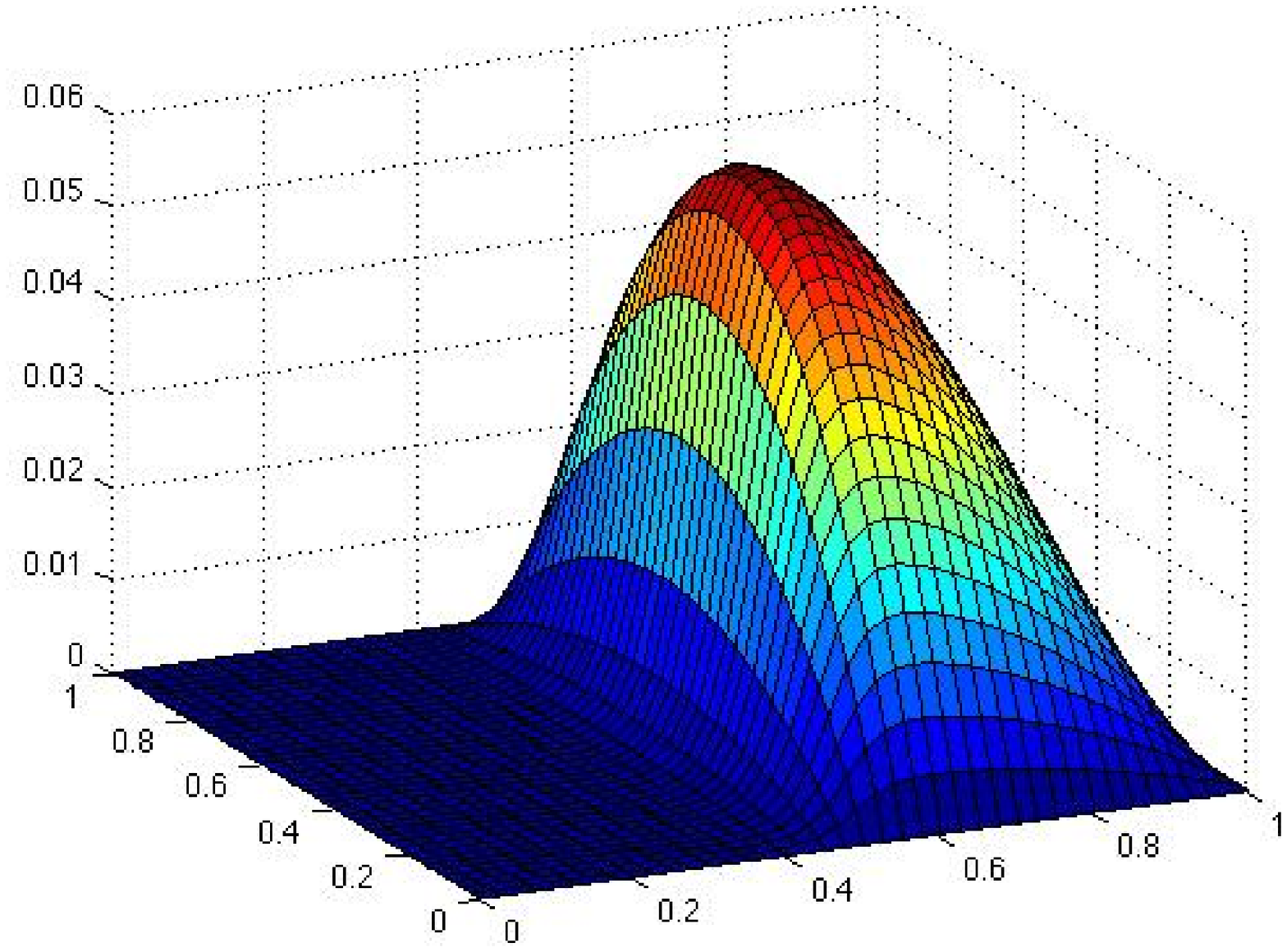}
\caption{The exact solution $u=\displaystyle \frac12 x(x-1)y(y-1)(1-tanh(10-20x))$.}
\label{exact_solution}
\end{center}
\end{figure}

Results are presented for uniformly refined meshes. In Tables
\ref{tableau_a1E-02} and \ref{tableau_a1E-04}, $N$ stands for the
number of mesh elements,
$\eta=\eta_{CF}+\eta_{NC}+\eta_{NC2}+\eta_{J}$, and
$Eff=\eta/||u-u_h||_{DG,h}$ is the effectivity index. $CV_{error}$
and $CV_{recov}$ are respectively the convergence order in
$\sqrt{1/N}$ of the error $||u-u_h||_{DG,h}$ and of
$||a^{-1/2}(Gu_h-a \nabla u)||$, from one line of the table to the
following. First of all, it can be observed for $\varepsilon=1E-02$
that the error $||u-u_h||_{DG,h}$ converges towards zero at order
one and that, in the same time, the superconvergence property of
$||a^{-1/2}(Gu_h-a \nabla u)||$ is observed.  Moreover, as expected
by Theorems \ref{tupperbound} and  \ref{tlowerbound2} when
$\bar\rho$ and $\bar\gamma$ are superconvergent terms, the proposed
estimator is reliable and efficient since the effectivity index
remains constant during the refinement process (around 2.00). For
$\varepsilon=1E-04$, the same kind of behaviour occurs. Let us
nevertheless note that the convergence rate $CV_{error}$ is
astonishment high. In fact, the contribution of the jump term
$\left(\sume h_e^{-1} \|\jump{u_h}\|^2_e\right)^{1/2}$ arising in
$||u-u_h||_{DG,h}$ is predominant  and goes fast towards zero, while
the contribution of $|||u-u_h|||$ is smaller but converges at order
one once the mesh is fine enough.
\begin{table}[ht]
\centering
\begin{tabular}{|c|c|c|c|c|c|c|}
\hline
&&&&&&\\
$N$ & $||u-u_h||_{DG,h}$  &$CV_{error}$ &$\eta$ & $Eff$ &  $||a^{-1/2}(Gu_h-a \nabla u)||$& $CV_{recov}$\\
&&&&&&\\
\hline
512 & 2.22E-02& &4.47E-02 & 2.01 & 4.29E-03 & \\
\hline
2048 & 1.13E-02 & 0.98 & 2.29E-02 &  2.02 &  1.61E-03 & 1.41 \\
\hline
8192 &5.58E-03 & 1.02 & 1.13E-02 &  2.02 &  4.87E-04 & 1.73 \\
\hline
32768 &2.77E-03 & 1.04 & 5.57E-03 & 2.01 & 1.37E-04 & 1.83 \\
\hline
131072  & 1.38E-03 & 1.01 & 2.77E-03 & 2.01 & 3.96E-05 & 1.79 \\
\hline
\end{tabular}
\caption{Homogeneous case, $\varepsilon=1E-02$, $\gamma=250$, $\gamma_{a,e}=\varepsilon$.\label{tableau_a1E-02}}
\end{table}

\begin{table}[ht]
\centering
\begin{tabular}{|c|c|c|c|c|c|c|}
\hline
&&&&&&\\
$N$ & $||u-u_h||_{DG,h}$  &$CV_{error}$ &$\eta$ & $Eff$ &  $||a^{-1/2}(Gu_h-a \nabla u)||$& $CV_{recov}$\\
&&&&&&\\
\hline
512 & 2.60E+00& & 2.72E+00 & 1.05 & 3.52E-02 &  \\
\hline
2048 & 1.46E+00 & 0.82& 8.28E-01&  1.06&  2.17E-02 & 0.69 \\
\hline
8192 & 5.38E-01& 1.44& 5.81E-01&  1.09&  7.97E-03 &1.44 \\
\hline
32768 &1.84E-01 &1.54&2.07E-01 & 1.12&2.40E-03&1.73\\
\hline
131072 & 6.83E-02 & 1.43 & 8.06E-02 & 1.18 & 7.14E-04 & 1.75 \\
\hline
\end{tabular}
\caption{Homogeneous case, $\varepsilon=1E-04$, $\gamma=250$, $\gamma_{a,e}=\varepsilon$.\label{tableau_a1E-04}}
\end{table}

\subsection{The singular case}
\label{singular} We consider here the domain $\Omega=\{ -1<x,y<1\}$,
which is decomposed into 4 sub-domains $\Omega_i$, $i=1,...,4$, with
$\Omega_1=(0,1) \times (0,1)$, $\Omega_2=(-1,0) \times (0,1)$,
$\Omega_3=(-1,0) \times (-1,0)$ and $\Omega_4=(0,1) \times (-1,0)$.
Like in subsection \ref{regular}, the reaction coefficient is
$\mu=1$, the velocity field $\beta=(1,0)^\top$, but the isotropic
diffusion tensor is this time no more homogeneous. It is defined by
$a_{| \Omega_i}=\varepsilon_i \, \mathcal{I}$, with
$\varepsilon_2=\varepsilon_4=1$ and
$\varepsilon_1=\varepsilon_3=C>1$ to be specified. \\
\\
Using the usual polar coordinates $(r,\theta)$ centered at $(0,0)$, the exact solution is chosen to be equal to the
singular function $u(x,y)= \eta(r) \, v(x,y)$ with $v(x,y)=r^{\alpha} \phi(\theta)$, where $\alpha
\in (0,1)$ and $\phi$ are chosen such that $v$ is harmonic on each
sub-domain $\Omega_i$, $i=1,..4$, and satisfies the jump conditions
:
$$
\left[ v \right]=0 \mbox{ and } \left[  a \nabla v . n \right] = 0
$$
on the interfaces. The function $\eta$ is a $\mathcal{C}^1[0,1]$ truncation function used to ensure homogeneous Dirichlet boudary condition on the boundary. Namely,
we chose :
$$
\eta(r)=
\left\{
\begin{array}{lll}
1 & \mbox{ for } 0 \leq r \leq 1/3, \\
(r-2/3)^2 \, (54r-9) & \mbox{ for } 1/3 \leq r \leq 2/3, \\
0 &  \mbox{ for } r \geq 2/3.\\
\end{array}
\right.
$$
It is easy to see (see for instance
\cite{CDN}) that $\alpha$ is the root of the transcendental equation
$$\tan \frac{\alpha \, \pi}{4}= \sqrt{\frac{1}{C}},
$$
and since $\alpha < 1$, this solution has a singular behavior around the
point $(0,0)$. Consequently, a local mesh-refinement strategy
is used, based on the estimator $\eta_T=\eta_{CF,T}+\eta_{NC,T}+\eta_{J,T}$
 and the marking process
 $$
 \eta_T>0.75 \max_{T'}\eta_{T'},$$
 with a standard refinement procedure associated with a limitation on the
 minimal angle. \\
\\
Like in \cite{creuse:10}, we choose first $C=5$. Figure
\ref{fig_dis_sing_a5_adap} shows some of the meshes obtained during
the local refinement process. Moreover, Table \ref{tableau_sing_a5}
displays the corresponding quantitative results.  It can be observed that the error
goes towards zero as theoretically expected, and that the effectivity index
always remains almost constant,  which is quite satisfactory and comparable
with results from \cite{cochez:07,ern:07} as well as those
of the previous test in subsection \ref{regular}. The superconvergence property of $||a^{-1/2}(Gu_h-a \nabla u)||$ is once again observed, and the mesh is automatically refined in the vicinity of
the singularity as well as in the zone of the mesh where the gradient of $\eta$ is the highest.
\\

\begin{figure}[htp]
\centering \setlength{\unitlength}{1mm}
\begin{tabular}{ccc}
\includegraphics*[width=4cm,height=4cm,angle=270]{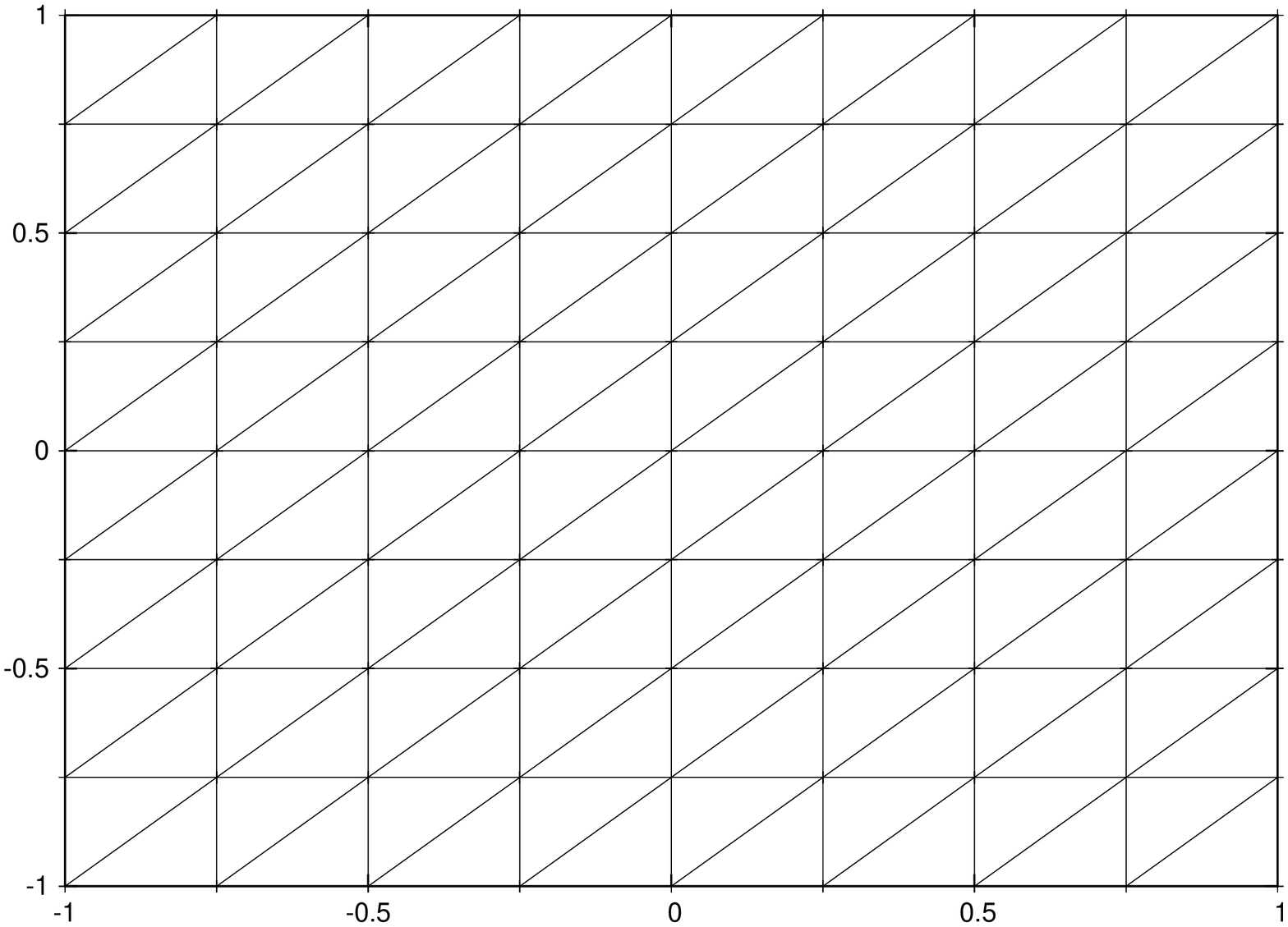} &
\includegraphics*[width=4cm,height=4cm,angle=270]{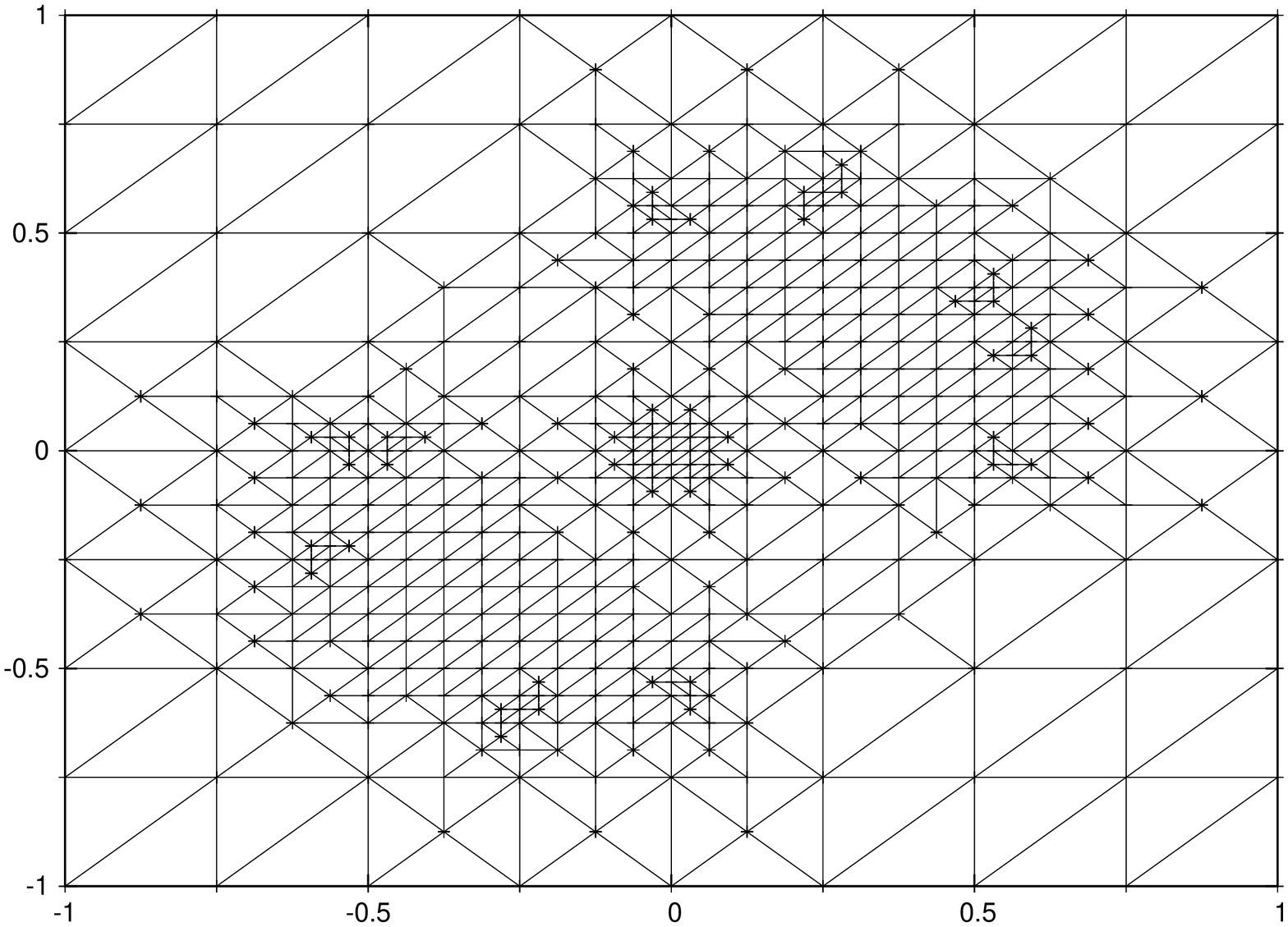} &
\includegraphics*[width=4cm,height=4cm,angle=270]{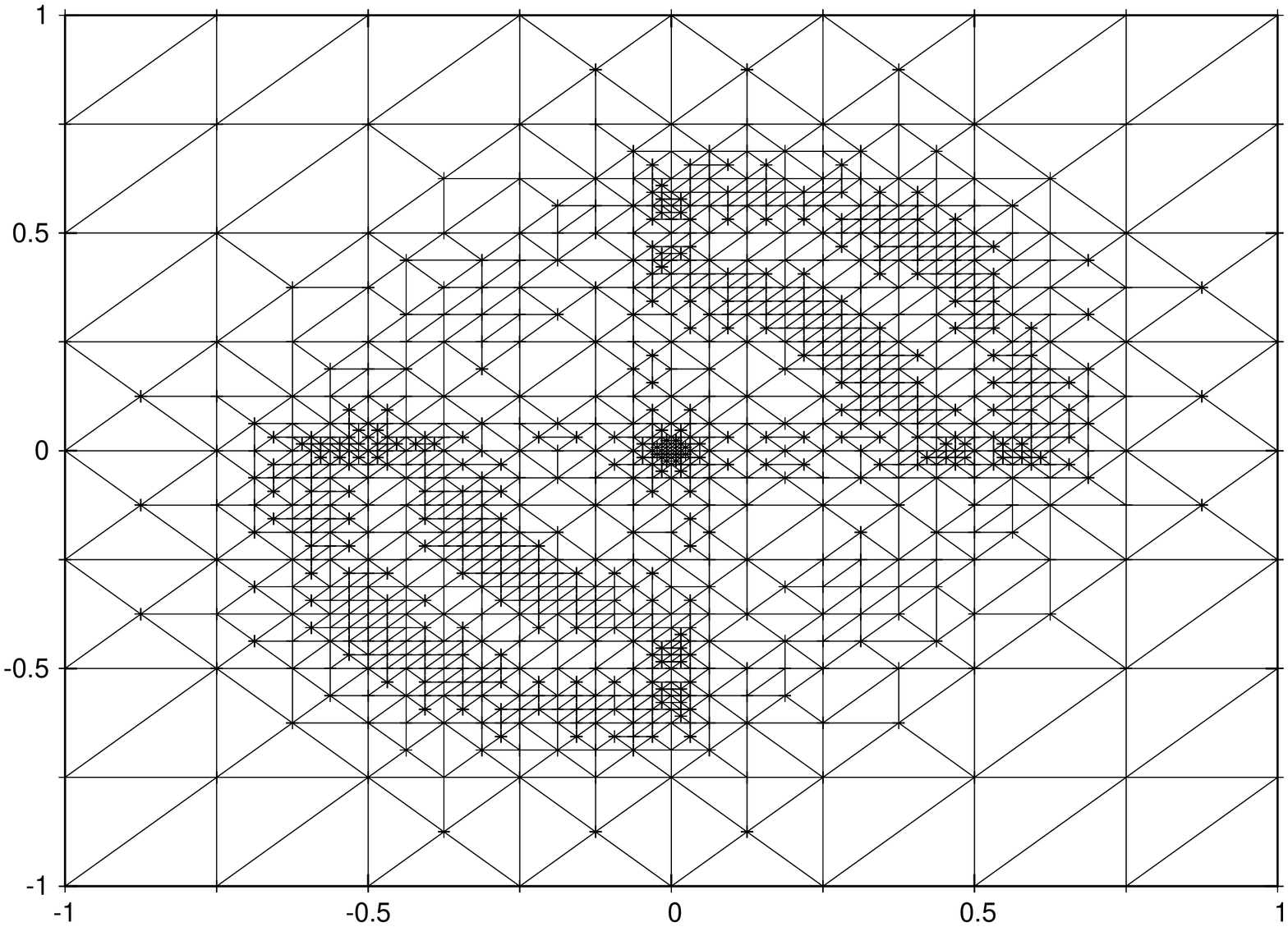}
\end{tabular}
\caption{Mesh levels 1, 3 and 5, singular solution for $C=5$.\label{fig_dis_sing_a5_adap}}
\end{figure}

\begin{table}[htp]
\centering
\begin{tabular}{|c|c|c|c|c|c|c|}
\hline
&&&&&&\\
$N$ & $||u-u_h||_{DG,h}$  &$CV_{error}$ &$\eta$ & $Eff$ &  $||a^{-1/2}(Gu_h-a \nabla u)||$& $CV_{recov}$\\
&&&&&&\\
\hline
128 &1.26E+00 & &3.09E+00 & 2.44  &1.14E+00& \\
\hline
468 &4.72E-01&1.51 & 1.25E+00&  2.65&  4.55E-01&1.41  \\
\hline
2016& 2.14E-01&1.08 & 5.74E-01&  2.68&  1.36E-01& 1.65 \\
\hline
9068 &1.01E-01 & 1.00&2.81E-01 & 2.76 & 4.68E-02 &1.41 \\
\hline
\end{tabular}
\caption{Homogeneous case, $C=5$, $\gamma=50$, $\gamma_{a,e}=1$.\label{tableau_sing_a5}}
\end{table}

Secondly, we consider a stronger singularity by choosing $C=100$.
Figure \ref{fig_dis_sing_a100_adap} shows some of the meshes
obtained during the local refinement process, and Table
\ref{tableau_sing_a100} displays the corresponding quantitative
results. The results are very similar to the ones obtained for
$C=5$. Like in \cite{creuse:10}, the refinement process is faster
around the interfaces (and the origin). The effectivity index
slightly increases while remaining constant during all the
refinement process, and the superconvergence property of
$||a^{-1/2}(Gu_h-a \nabla u)||$ is observed.

\begin{figure}[htp]
\centering \setlength{\unitlength}{1mm}
\begin{tabular}{ccc}
\includegraphics*[width=4cm,height=4cm,angle=270]{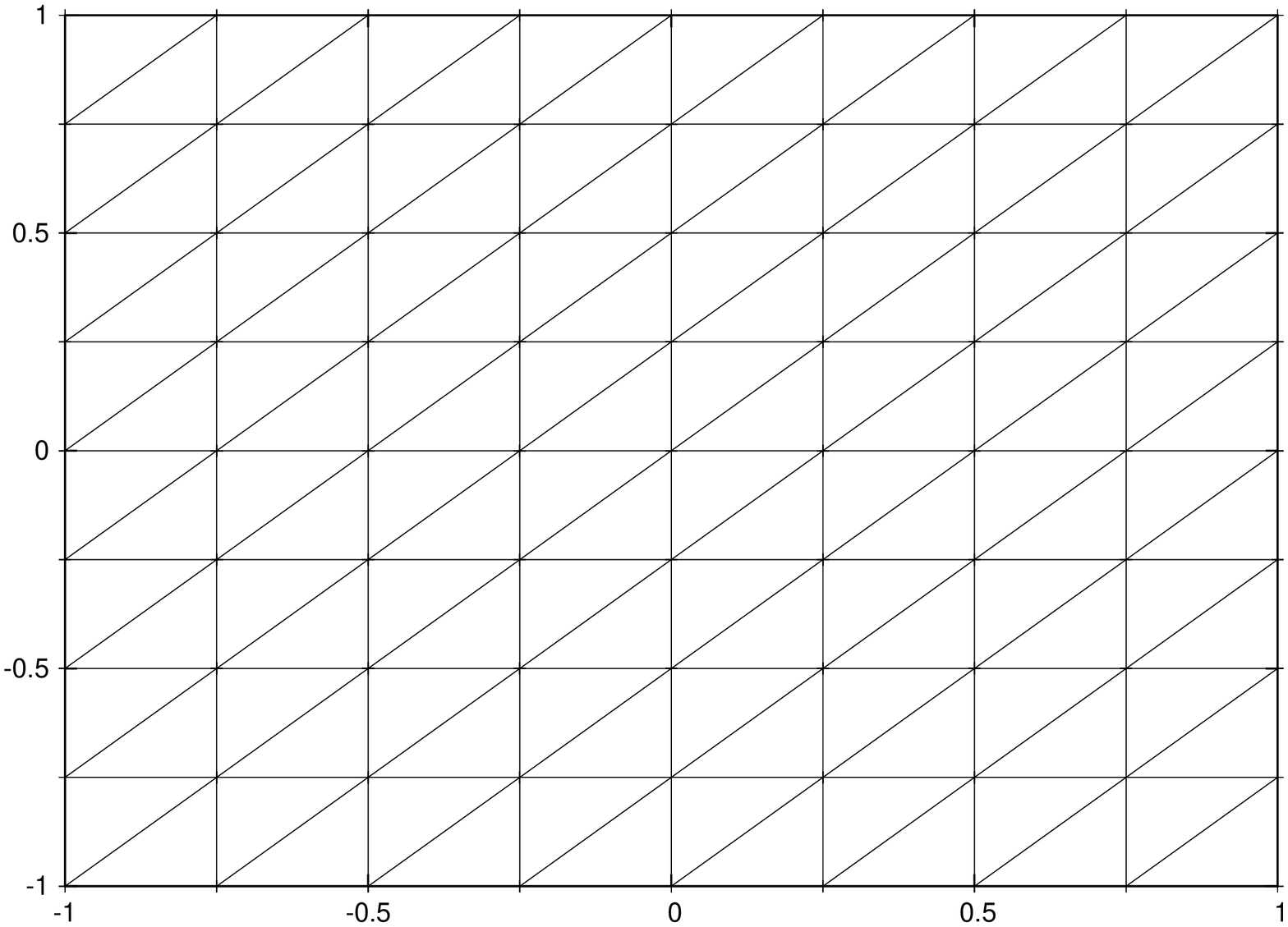} &
\includegraphics*[width=4cm,height=4cm,angle=270]{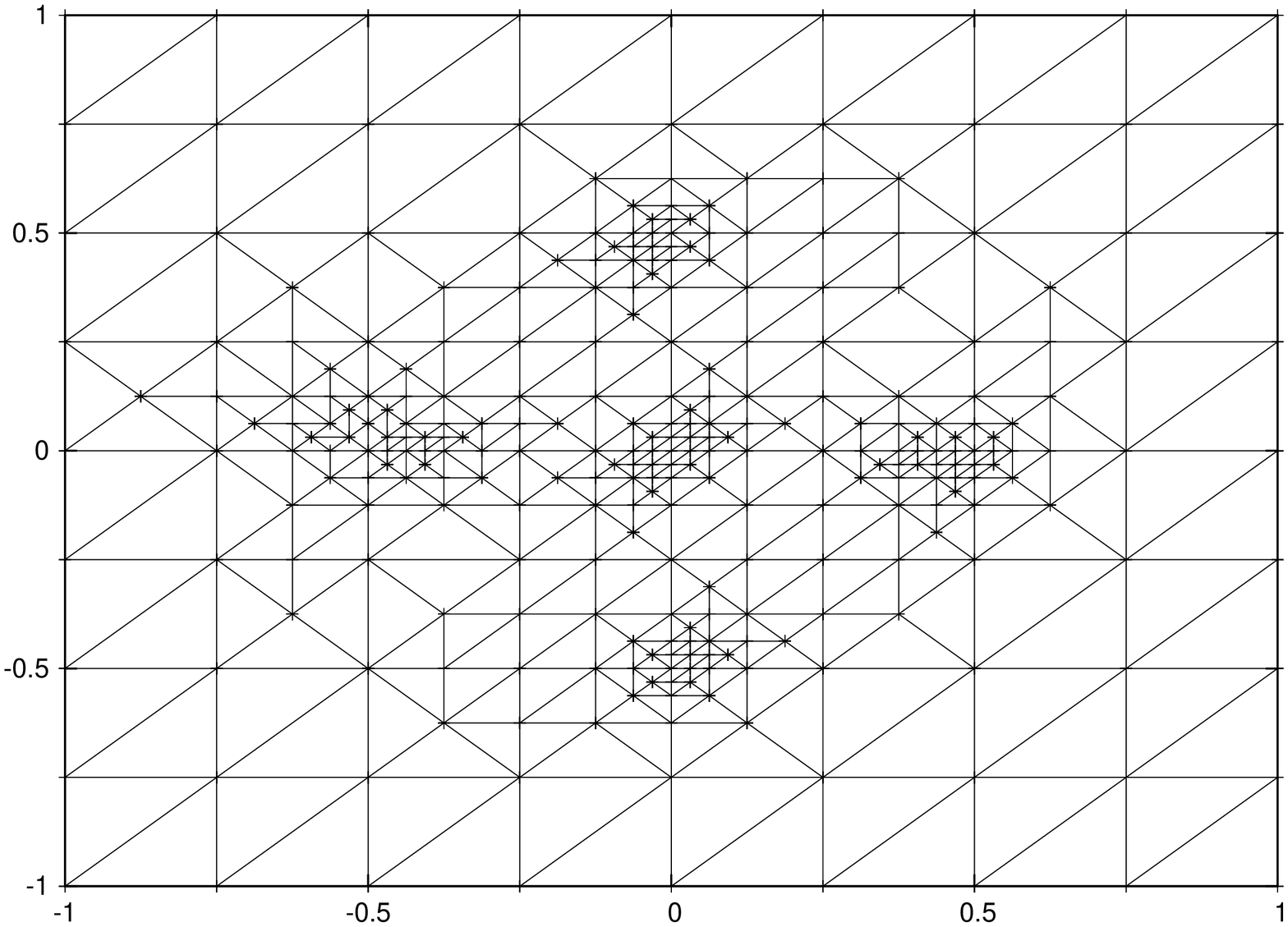} &
\includegraphics*[width=4cm,height=4cm,angle=270]{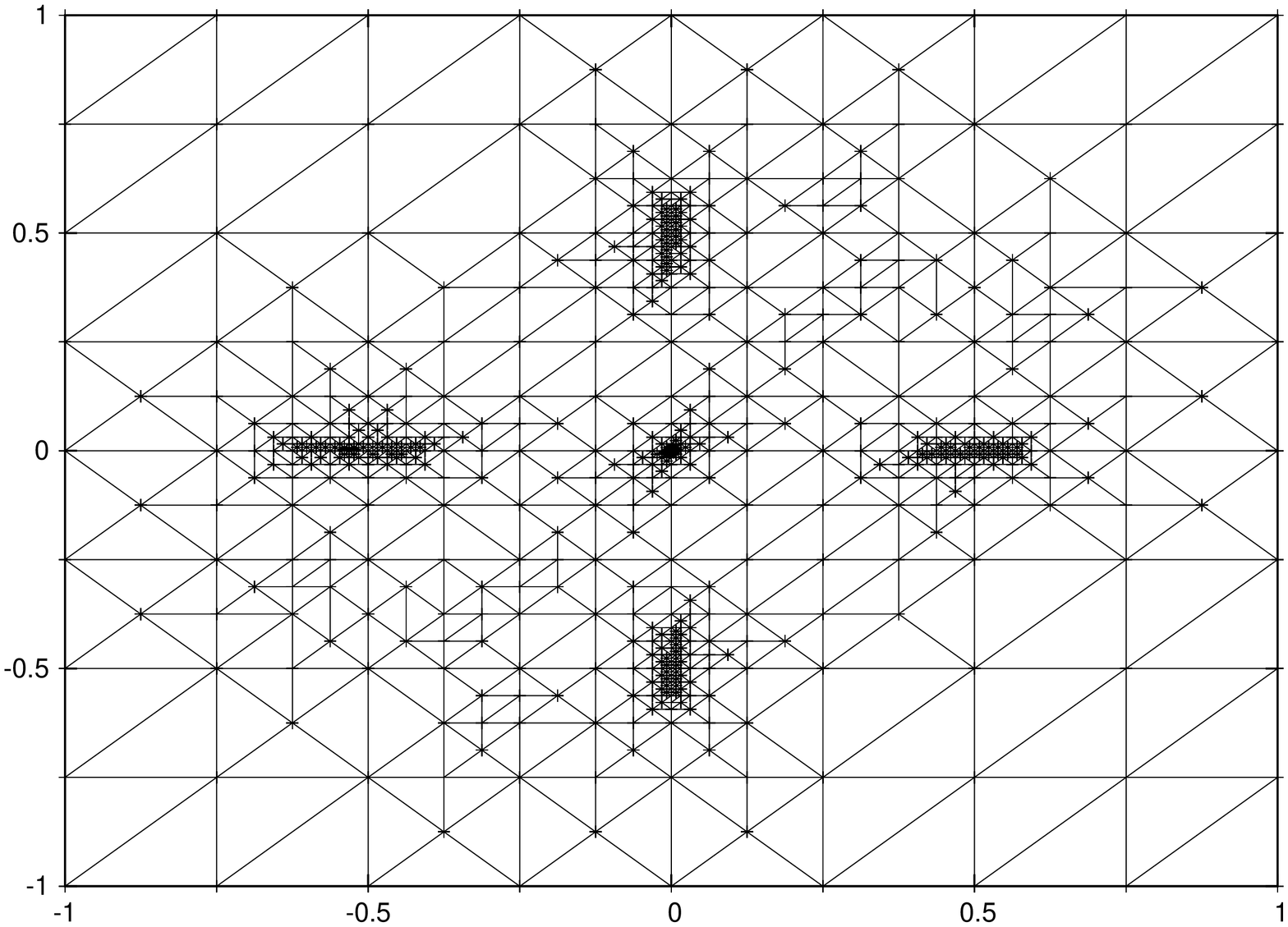}
\end{tabular}
\caption{Mesh levels 1, 3 and 15, singular solution for $C=100$.\label{fig_dis_sing_a100_adap}}
\end{figure}

\begin{table}[ht]
\centering
\begin{tabular}{|c|c|c|c|c|c|c|}
\hline
&&&&&&\\
$N$ & $||u-u_h||_{DG,h}$  &$CV_{error}$ &$\eta$ & $Eff$ &  $||a^{-1/2}(Gu_h-a \nabla u)||$& $CV_{recov}$\\
&&&&&&\\
\hline
128 &1.56E+00 & &6.34E+00 & 4.05  &2.92E+00& \\
\hline
504 &7.25E-01&1.12 & 3.65E+00&  5.03&  7.99E-01&1.89  \\
\hline
2208& 4.48E-01&0.65& 1.95E+00&  4.35&  3.81E-01& 1.00 \\
\hline
7716 & 2.44E-01 &0.97 & 1.09E+00 & 4.46 & 1.60E-01  &1.38 \\
\hline
\end{tabular}
\caption{Homogeneous case, $C=100$,  $\gamma=500$, $\gamma_{a,e}=1$.\label{tableau_sing_a100}}
\end{table}
\subsection{The  boundary-layer case}
We now consider the domain $\Omega=\{ 0<x,y<1\}$, with the reaction coefficient $\mu=0$ and the velocity
field $\beta=(1,0)^\top$.  The homogeneous isotropic diffusion tensor is defined by $a=\varepsilon \, \mathcal{I}$, with $\varepsilon= 1E-02$. The source term $f$ is chosen accordingly so that $u=10y(1-y)x(e^{-x}-e^{-1+\frac{x-1}{\alpha}})$ is the exact solution (see Figure \ref{exact_solution_cl}), with $\alpha=1E-03$ in order to generate a strong boundary layer.  Here, the same refinement procedure than in section \ref{singular}  is used. \\

\begin{figure}[htp]
\begin{center}
\includegraphics[scale=0.50]{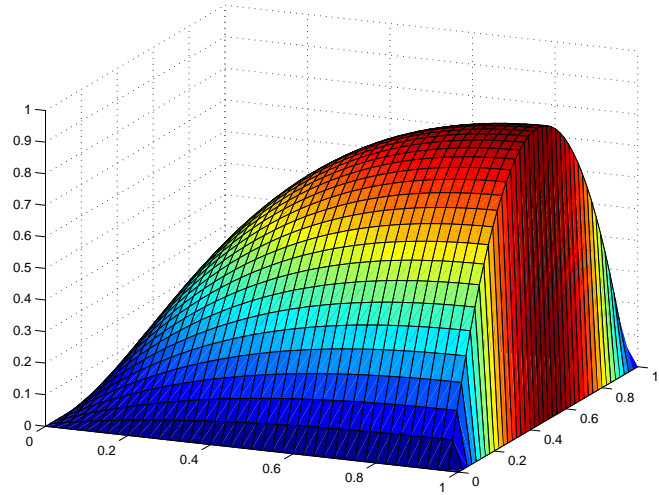}
\caption{The exact solution $u=10y(1-y)x(e^{-x}-e^{-1+\frac{x-1}{\alpha}})$, with $\alpha=1E-03$.}
\label{exact_solution_cl}
\end{center}
\end{figure}

Figure \ref{meshes_singconv} shows some of the meshes obtained during
the local refinement process. Moreover, Table \ref{tableau_sing_convaE-02}
displays the corresponding quantitative results.  Provided that the boundary layer mesh resolution is
sufficient, the same behaviours than in the previous tests  can also be observed : the error
goes towards zero as theoretically expected, the effectivity index
remains almost constant, the superconvergence property of $||a^{-1/2}(Gu_h-a \nabla u)||$ occurs, and  the mesh is automatically refined around the boundary layer.

\begin{figure}[htp]
\centering \setlength{\unitlength}{1mm}
\begin{tabular}{ccc}
\includegraphics*[width=4cm,height=4cm,angle=270]{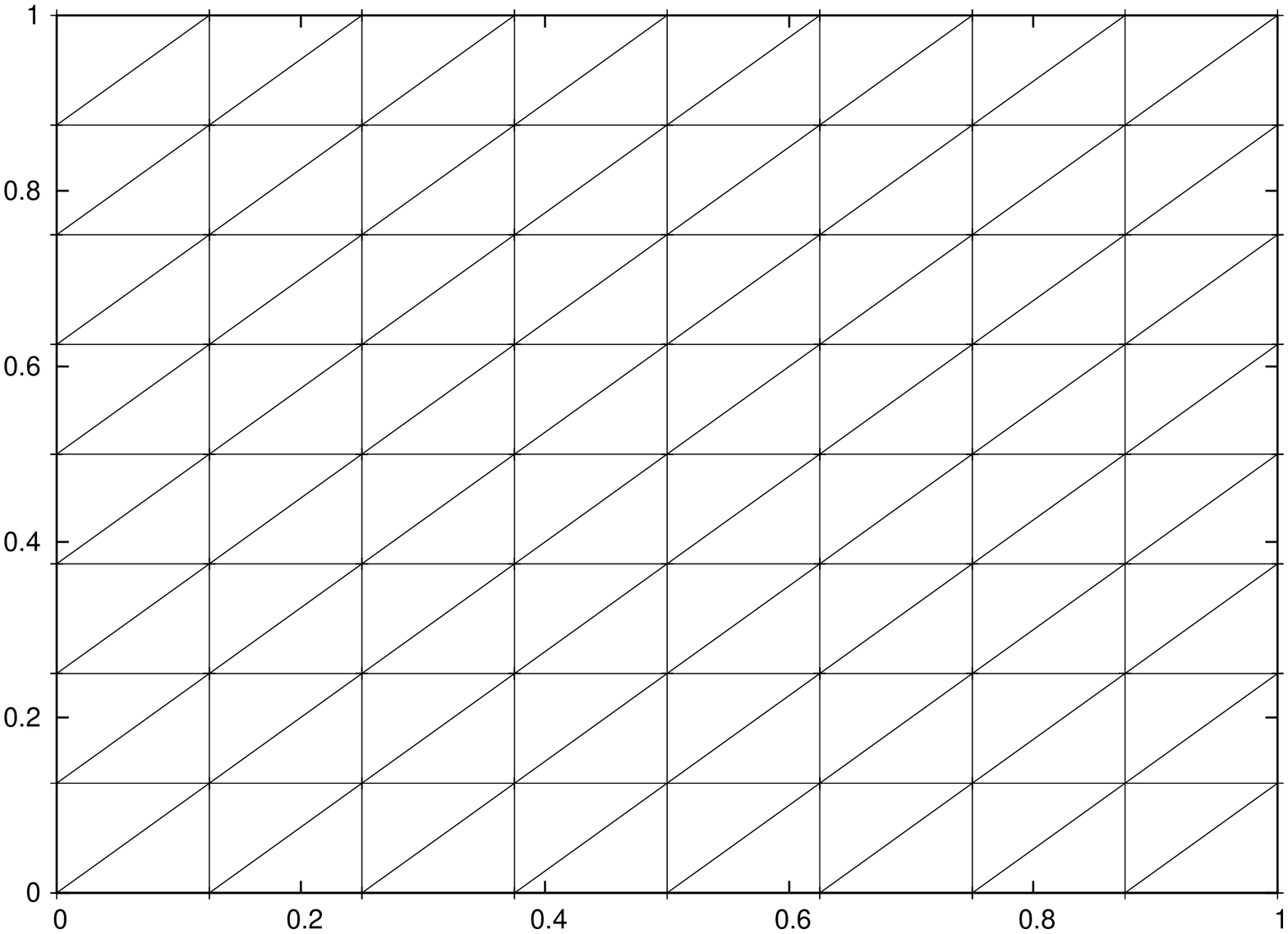} &
\includegraphics*[width=4cm,height=4cm,angle=270]{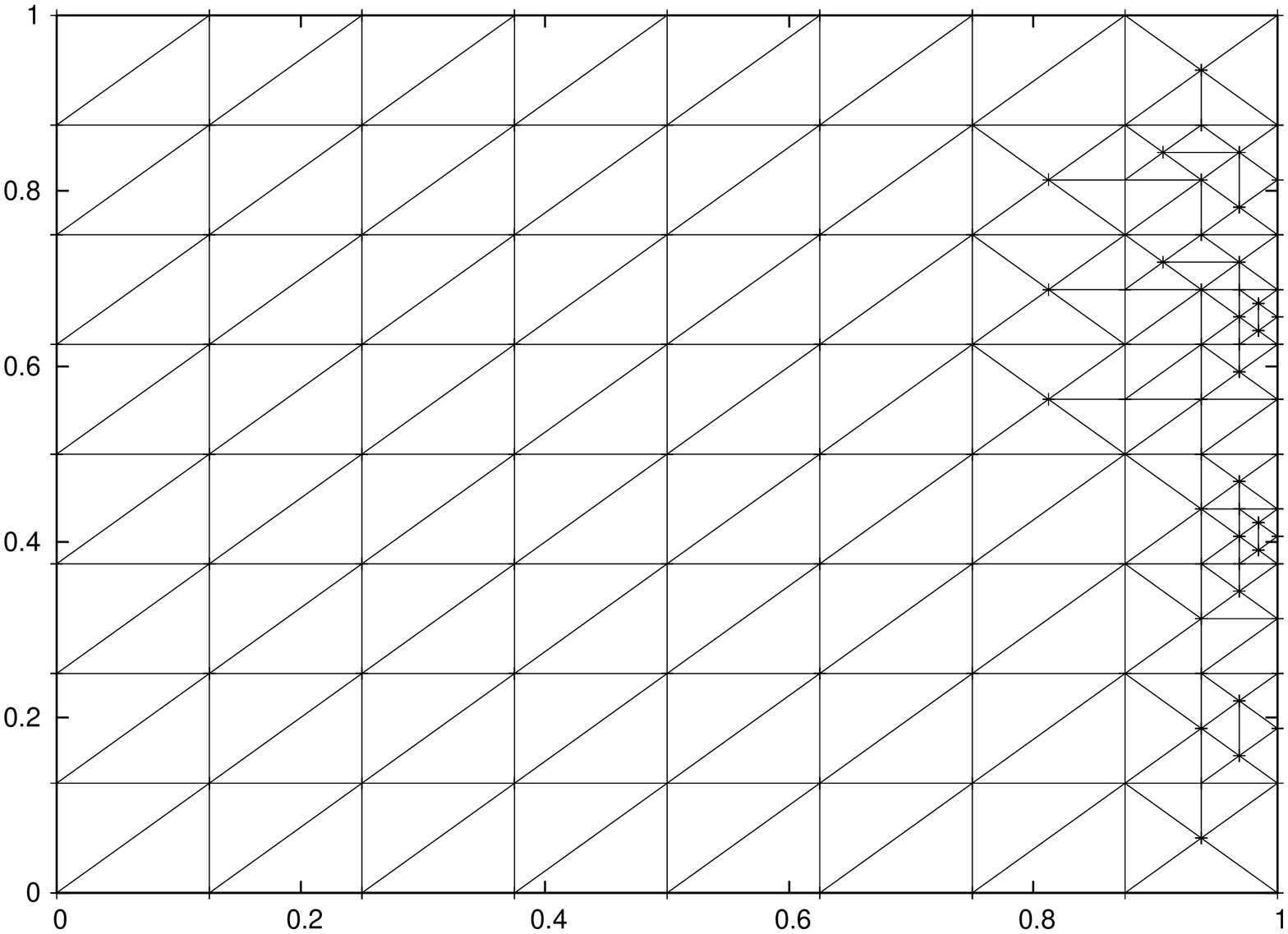} &
\includegraphics*[width=4cm,height=4cm,angle=270]{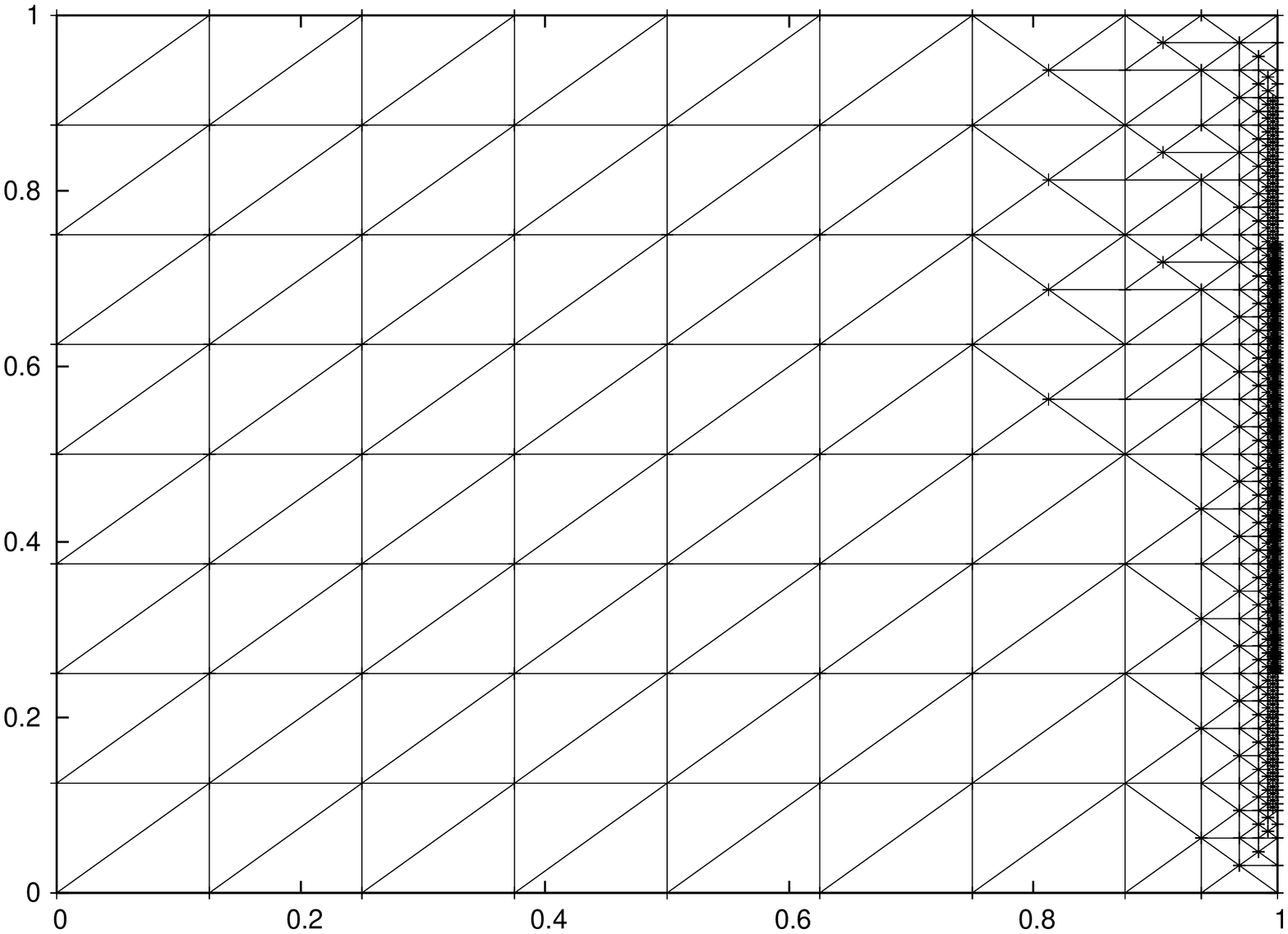}
\end{tabular}
\caption{Mesh levels 1, 5 and 16, non homogeneous boundary-layer case.\label{meshes_singconv}}
\end{figure}

\begin{table}[ht]
\centering
\begin{tabular}{|c|c|c|c|c|c|c|}
\hline
&&&&&&\\
$N$ & $||u-u_h||_{DG,h}$  &$CV_{error}$ &$\eta$ & $Eff$ &  $||a^{-1/2}(Gu_h-a \nabla u)||$& $CV_{recov}$\\
&&&&&&\\
\hline
 156 & 5.90E+00& & 9.11E+00& 1.54&  1.46E+00 & \\
 \hline
766 & 1.32E+00 & 1.88 & 2.52E+00 & 1.90 & 8.84E-01 & 0.63 \\
\hline
3033 & 5.76E-01 & 1.21 & 1.24E+00 & 2.16 & 3.67E-01 & 1.27 \\
\hline
12881 & 2.53E-01 & 1.13& 5.34E-01 &  2.10 & 1.19E-01 & 1.58\\
\hline
50496 & 1.26E-01 & 1.02 &2.66E-01 & 2.11 & 5.02E-02 &1.26\\
\hline
\end{tabular}
\caption{Boundary-layer case, $\varepsilon=1E-02$, $\gamma=250$, $\gamma_{a,e}=\varepsilon$. \label{tableau_sing_convaE-02}}
\end{table}

\clearpage

    \protect\bibliographystyle{abbrv}
    \protect\bibliographystyle{alpha}
    \bibliography{E:/documents/Serge/Desktop/Biblio/kunert,E:/documents/Serge/Desktop/Biblio/est,E:/documents/Serge/Desktop/Biblio/femaj,E:/documents/Serge/Desktop/Biblio/mgnet,E:/documents/Serge/Desktop/Biblio/dg,E:/documents/Serge/Desktop/Biblio/bib,E:/documents/Serge/Desktop/Biblio/maxwell,E:/documents/Serge/Desktop/Biblio/bibmix,E:/documents/Serge/Desktop/Biblio/cochez,E:/documents/Serge/Desktop/Biblio/soualem,E:/documents/Serge/Desktop/Biblio/nic}


\end{document}